\documentclass[lettersize,journal]{IEEEtran}
\usepackage{cite}
\usepackage{amssymb,amsfonts}
\usepackage{algorithm}
\usepackage{algorithmic}
\usepackage{graphicx}
\usepackage{textcomp}
\usepackage{xcolor}

\newcounter{counter}

\newtheorem{theorem}{Theorem}
\newtheorem{definition}{Definition}
\newtheorem{proposition}{Proposition}
\usepackage{booktabs, makecell}
\setcellgapes{3pt}
\usepackage{siunitx} 
\usepackage{multirow}
\usepackage{multicol}
\usepackage[caption=false,font=normalsize,labelfont=sf,textfont=sf]{subfig}
\usepackage{fancyhdr}
\usepackage{stfloats}
\usepackage{url}
\usepackage{verbatim}
\usepackage{amsmath}
\usepackage{amssymb}
\usepackage{bm}
\usepackage{textcomp}

\makeatletter
\newenvironment{breakablealgorithm}
  {
   \begin{center}
     \refstepcounter{algorithm}
     \hrule height.8pt depth0pt \kern2pt
     \renewcommand{\caption}[2][\relax]{
       {\raggedright\textbf{\ALG@name~\thealgorithm} ##2\par}%
       \ifx\relax##1\relax 
         \addcontentsline{loa}{algorithm}{\protect\numberline{\thealgorithm}##2}%
       \else 
         \addcontentsline{loa}{algorithm}{\protect\numberline{\thealgorithm}##1}%
       \fi
       \kern2pt\hrule\kern2pt
     }
  }{
     \kern2pt\hrule\relax
   \end{center}
  }
\makeatother
\newcommand{\Call}[2]{\textsc{#1}(#2)}
\begin{document}

\title{A Dynamic Relaxation Framework for Global Solution of ACOPF}

\author{Yu-Yang Tang, Liang Chen, Sheng-Jie Chen, Yu-Hong Dai,\\ Bo Zhou,~\IEEEmembership{Member,~IEEE} and Xiaomeng Ai,~\IEEEmembership{Member,~IEEE,}
\thanks{This research was supported by the National Key R\&D Program of China (No.
2022YFB2403400) and the Chinese NSF grants (No. 12201620). (\emph{Corresponding author: Liang Chen})}
\thanks{Yu-Yang Tang, Liang Chen, Sheng-Jie Chen and Yu-Hong Dai are with the Institute of Computational Mathematics and Scientific/Engineering Computing, Academy of Mathematics and Systems Science, Chinese Academy of Sciences, Beijing 100190, China (e-mail: tangyuyang@lsec.cc.ac.cn, chenliang@lsec.cc.ac.cn, shengjie\_chen@lsec.cc.ac.cn, dyh@lsec.cc.ac.cn).}
\thanks{Bo Zhou is with the Department of Industrial and Operations Engineering, University of Michigan, Ann Arbor, MI 48109 USA (e-mail: bozum@umich.edu).}
\thanks{Xiaomeng Ai is with the State Key Laboratory of Advanced Electromagnetic Engineering and Technology, School of Electrical and Electronic Engineering,
Huazhong University of Science and Technology, Wuhan 430074, China (e-mail: xiaomengai@hust.edu).}
}

\markboth{IEEE Transactions on Power Systems,~Vol.~xx, No.~x, xx~2025}
{Tang \MakeLowercase{\textit{et al.}}: A Dynamic Relaxation Framework for Global Solution of ACOPF}


\maketitle

\begin{abstract}
Solving the Alternating Current Optimal Power Flow (AC OPF) problem to global optimality remains challenging due to its nonconvex quadratic constraints. In this paper, we present a unified framework that combines static piecewise relaxations with dynamic cut-generation mechanism to systematically tighten the classic Second-Order Cone Programming (SOCP) relaxation to arbitrarily small conic violation, thus enabling the recovery of globally optimal solutions. Two static formulations, Pyramidal Relaxation (PR) and Quasi‐Pyramidal Relaxation (QPR), are introduced to tighten each branch‐flow second-order cone via a finite union of wedges, providing controllable accuracy. Their dynamic counterparts, Dynamic PR (DPR) and Dynamic QPR (DQPR), embed on‐the‐fly cut generation within a branch‐and‐cut solver to improve scalability. Convergence is further accelerated through warm starts and a lightweight local‐search post‐processing. Extensive experiments on benchmarks demonstrate effective elimination of conic violations and flexible trade‐offs between solution accuracy and runtime. Practical guidelines are derived for selecting appropriate variants based on network size and accuracy requirements. 
\end{abstract}

\begin{IEEEkeywords}
Alternating Current Optimal Power Flow (AC OPF), Second-Order Cone Surface Programming (SOCSP), Pyramidal Relaxation (PR), Quasi-Pyramidal Relaxation (QPR), Mixed-Integer Linear Programming (MILP), Mixed-Integer Second-Order Cone Programming (MISOCP),  dynamic relaxtion framework.
\end{IEEEkeywords}

\vspace{-2ex}
\section{Introduction}
%
%
%
%

\IEEEPARstart{T}{he} Optimal Power Flow (OPF) problem, first formulated in the 1960s, plays a foundational role in modern power system operations \cite{carpentier1962contribution}.
At its core, the OPF problem aims to determine the optimal settings for various controllable resources in the power grid to achieve objectives such as minimizing generation costs, while adhering to physical and operational constraints. This problem is crucial for enhancing the efficiency, reliability, and economic viability of power systems \cite{stott2012optimal}. As a cornerstone of electrical engineering, OPF facilitates critical functions ranging from security-constrained optimizations \cite{wu2018robust} to voltage stability assessments \cite{mohamed2019voltage}.

Despite its importance, the OPF problem, especially in its Alternating Current (AC) form, is inherently difficult to solve due to its non-linear and complex nature. Moreover, it has been proven to be NP-hard \cite{lavaei2011zero}, posing fundamental challenges to global optimization. While several modeling paradigms exist, including the widely used Bus Injection Model (BIM), this work adopts the Branch Flow Model (BFM)~\cite{baran1989optimal1,farivar2013branch}, which explicitly represents branch power flows and currents. For initialization, we also adopt the current–voltage (IV) formulation in rectangular coordinates~\cite{o2012iv} to generate high-quality warm-start solutions.


Throughout the long history of OPF research, numerous optimization techniques have been developed to address its nonconvexity. Interior Point Methods (IPMs)~\cite{wu1994direct} are widely used in practice, but cannot guarantee global optimality due to the NP-hard nature of AC OPF~\cite{lavaei2011zero}. To overcome this, convex relaxations have attracted significant attention. The Semidefinite Programming (SDP) relaxation was pioneered by Bai~{\it et al.}~\cite{bai2008semidefinite} and popularized by Lavaei and Low~\cite{lavaei2011zero}, who established sufficient conditions for exactness. However, SDP often becomes computationally prohibitive as network size grows, and its exactness can fail under tight line constraints~\cite{lesieutre2011examining}. Refinements such as moment-based hierarchy relaxations~\cite{molzahn2014sparsity} and global branch-and-bound methods~\cite{phan2012lagrangian} offer tighter bounds at higher cost.
As a more tractable alternative, Second-Order Cone Programming (SOCP) relaxations have been proposed. Jabr~\cite{jabr2006radial} introduced a BIM-based SOCP for radial networks, while Farivar and Low~\cite{farivar2013branch} extended this to mesh networks via the BFM. Yet, exact conditions such as load oversatisfaction remain unrealistic in general networks. Recent efforts aim to improve both scalability and tightness of the relaxtation, including Quadratic Relaxations (QC)~\cite{coffrin2015qc}, strong SOCP~\cite{kocuk2016strong}, and learning-based methods~\cite{pan2020deepopf}. For broader reviews, see~\cite{molzahn2019survey}.

Linear approximations have long been used to simplify the AC OPF problem for improved tractability. Early approaches such as DC OPF~\cite{stott2009dc} are widely adopted in practice, but may be inaccurate under system stress. More refined models based on local linearization—using polar, rectangular, or IV formulations~\cite{yang2017linearized, dhople2015linear, castillo2015successive}—offer better fidelity but depend heavily on the linearization point. To improve robustness, piecewise linear approximations~\cite{zhou2020pyramidal, akbari2016linear} yield MILP formulations, though their combinatorial complexity limits scalability. Meanwhile, convex linear relaxations~\cite{coffrin2016network, coffrin2014linear} provide dual bounds but often sacrifice accuracy due to simplifying assumptions.

As discussed above, the sufficient conditions required for exact convexification of OPF are rarely satisfied in real-world scenarios. In particular, SOCP relaxations often yield solutions that lie strictly within the interior of the cone, rather than on its boundary, resulting in overestimated power losses and other inaccuracies~\cite{zhou2020pyramidal}. To address this inexactness, several refinement schemes have been explored. Liu~{\it et al.} proposed an outer-approximation method based on McCormick envelopes~\cite{liu2018multitree}, but the resulting bounds are often loose, and the formulation scales poorly. Zhou~{\it et al.} introduced piecewise linear constraints to approximate the SOC surface~\cite{zhou2020pyramidal}, improving accuracy but risking infeasibility when the approximated region fails to intersect the true cone.

Building on our prior work~\cite{tang2024dynamic}, this paper proposes a dynamic relaxation framework that combines the tightness and robustness of piecewise relaxations with a branch-and-cut mechanism equipped with on-the-fly cut generation. This approach improves solution accuracy and scalability, while maintaining feasibility and offering tighter dual bounds.

Our main contributions are as follows:
\begin{itemize}
    \item We introduce the second-order cone surface programming (SOCSP) formulation for AC OPF and define associated relative and absolute conic error metrics. Based on this, we propose two static relaxations—Pyramidal Relaxation (PR) and Quasi-Pyramidal Relaxation (QPR)—which partition the branch-flow cone surface into finitely many wedges and are provably asymptotically exact.
    \item We develop a branch-and-cut–based dynamic relaxation framework that starts from minimal models and generates only the necessary cuts during solving. Warm-start and post-processing techniques are integrated to further enhance performance.
    \item We perform extensive experiments on eight PGLib–OPF instances to demonstrate the effectiveness of proposed relaxations and dynamic methods. Based on the results, we offer practical guidance on method selection with respect to network size and desired conic accuracy.
\end{itemize}

The remainder of this paper is organized as follows. Section~\ref{sec:pf} presents the AC OPF formulation, introduces the SOCSP model, and defines conic error metrics. Section~\ref{sec:pr} describes the static formulations (PA, PR, QPR). Section~\ref{sec:drf} details the dynamic relaxation framework. Numerical results are reported in Section~\ref{sec:nr}, and Section~\ref{sec:concl} concludes the paper. 

\emph{This is the full version of a paper submitted to IEEE Transactions on Power Systems. It includes all proofs and algorithmic pseudocode omitted in the main submission due to space limits.}
\vspace{-2ex}
\section{Problem Formulation}\label{sec:pf}
This section presents the specific AC OPF formulation used in our study, based on the relaxed branch flow model 
\cite{farivar2013branch}, derives its SOCP relaxation and defines metrics to quantify any resulting inexactness.\vspace{-2.5ex}
\subsection{Notation}
Let \( \mathcal{N} \), \( \mathcal{L} \), and \( \mathcal{G} \) denote the sets of buses, branches, and generators, respectively. A branch \( l \in \mathcal{L} \) is directed from bus \( i \) to \( j \), denoted \( l = i \to j \). The following variables are used in the optimization model: \( p_g \), \( q_g \) are the active and reactive power outputs of generator \( g \in \mathcal{G} \); \( V_i \) is the voltage magnitude at bus \( i \in \mathcal{N} \); \( P_{ij}, Q_{ij}, S_{ij}, I_{ij} \) represent the active power, reactive power, apparent power, and current magnitude on branch \( i \to j \), respectively. For notational simplicity, we may also write \( P_l = P_{ij} \), \( Q_l = Q_{ij} \), and similarly for \( S_l \), \( I_l \) when branch \( l \) is clear from context.

The model parameters include the branch resistance \( r_{ij} \), reactance \( x_{ij} \), the demand \( p_{i}^d, q_{i}^d \) at each bus \( i \), and generator limits \( p_g^{\min}, p_g^{\max}, q_g^{\min}, q_g^{\max} \). Voltage and line flow limits are denoted by \( V_i^{\min}, V_i^{\max} \), and \( S_{ij}^{\max} \), respectively.
\vspace{-2.5ex}
\subsection{AC OPF Formulation}
The relaxed branch flow model (also known as OPF-ar~\cite{farivar2013branch}) is given by:
\begin{subequations}
    \begin{equation}
        \min\quad\sum_{g\in \mathcal{G}} F_{g}\left(p_{g}\right) \label{obj}
    \end{equation}
    \text{s.t.}
    \begin{equation}
        \sum_{g\in \mathcal{G}_i} p_{g}-p_{i}^d=\sum_{j:i\to j \in \mathcal{L}}P_{ij} - \sum_{k:k\to i \in \mathcal{L}}(P_{ki} - r_{ki}I_{ki}^2),\quad \forall i \in \mathcal{N} \label{powbal_p}
    \end{equation}
    \begin{equation}
         \sum_{g\in \mathcal{G}_i} q_{g}-q_{i}^d=\sum_{j:i\to j \in \mathcal{L}}Q_{ij} - \sum_{k:k\to i \in \mathcal{L}}(Q_{ki} - x_{ki}I_{ki}^2),\quad \forall i \in \mathcal{N} \label{powbal_q}
    \end{equation}
    \begin{equation}
        V_{i}^2-V_{j}^2=2 r_l P_{l}+2 x_l Q_{l}-\left(r_l^2+x_l^2\right) I_l^2,\quad\forall l=i\to j \in \mathcal{L}\label{branch_v}
    \end{equation}
    \begin{equation}
        S_{l}^2=P_{l}^2+Q_{l}^2=I_{l}^2 V_{i}^2,\quad\forall l=i\to j\in\mathcal{L} \label{socp}
    \end{equation}
    \begin{equation}
        p_{g}^{\min } \leq p_{g} \leq p_{g}^{\max }, q_{g}^{\min } \leq q_{g} \leq q_{g}^{\max },\quad \forall g\in \mathcal{G}\label{gen_bounds}
    \end{equation}
    \begin{equation}
        V_i^{\min } \leq V_i \leq V_i^{\max },\quad \forall i\in \mathcal{N}\label{vol_bounds}
    \end{equation}
    \begin{equation}
        0\leq S_{l}\leq S_{l}^{\max},\quad \forall l\in \mathcal{L}\label{flow_bounds}
    \end{equation}
    \begin{equation}
        \tan\theta_l^{\min}(V_i^2 - r_lP_l- x_lQ_l)\le x_lP_l-r_lQ_l, \quad\forall l=i\to j \in \mathcal{L}\label{ang_up}
    \end{equation}
      \begin{equation}
        x_lP_l-r_lQ_l\le \tan\theta_l^{\max}(V_i^2 - r_lP_l- x_lQ_l), \quad\forall l=i\to j \in \mathcal{L}.\label{ang_low}
    \end{equation}
    \label{formulation}
\end{subequations}\noindent
Here, $F_{g}(p_{g})$ denotes the generation cost of unit $g$, assumed linear for simplicity in this paper.  
The objective function \eqref{obj} minimizes the total generation cost. Equations \eqref{powbal_p} and \eqref{powbal_q} enforce active/reactive power balance at each bus. Equations \eqref{branch_v} are derived from the relationship between bus voltages and branch power flows, while~\eqref{socp} enforces the power–current–voltage relationship. Constraints \eqref{gen_bounds}-\eqref{flow_bounds} impose generator, voltage, and line‐flow limits, where $V_i$ and $S_{l}$ are nonnegative. \eqref{ang_up}–\eqref{ang_low} bound the phase‐angle difference within $(-\tfrac{\pi}{2},\tfrac{\pi}{2})$ by linearizing its tangent.

Formulation \eqref{formulation} omits the angle recovery condition of the original branch flow model; it can be reinstated by introducing virtual phase-shifter to the network as in \cite{farivar2013branch2}.  Since this paper focuses exclusively on the inexactness introduced by SOCP relaxations, angle recovery is excluded. Likewise, network shunt elements are omitted for notational brevity, though included in all numerical experiments.

Following \cite{farivar2013branch}, we substitute $W_i = V_i^2$ and $\Phi_l = I_l^2$ which linearize constraints \eqref{powbal_p}, \eqref{powbal_q}, \eqref{branch_v}, \eqref{ang_low}, and \eqref{ang_up}, leaving the nonlinear equality \eqref{socp} as the only remaining source of nonconvexity. The classical SOCP relaxation replaces~\eqref{socp} with the conic inequality: 
\begin{equation}
    P_l^2+Q^2_l \le \Phi_l W_i,
    \label{socp_proto}
\end{equation}
yielding a convex model known as OPF-cr \cite{farivar2013branch}. To handle the apparent power variable \(S_l\) more explicitly, we adopt the 3D SOC reformulation:
\begin{subequations}
    \begin{align}
    P_{l}^2 + Q_{l}^2 &= S_{l}^2, \label{cone1}\\
    S_{l}^2+\left(\frac{W_{i}-\Phi_{l}}{2}\right)^2&=\left(\frac{W_{i}+\Phi_{l}}{2}\right)^2,
    \label{cone2}
\end{align}\label{socp_new}
\end{subequations}
Each of these defines an \emph{SOC surface}; relaxing the equalities with ``\(\le\)"-inequalities yields a standard SOCP. 
A key property of SOCP relaxations is {\it exactness}: the relaxed solution lies on the SOC surfaces. Inexactness occurs when the solution lies strictly inside the cones, leading to infeasibility or suboptimality. Existing exactness conditions often rely on load‐oversatisfaction assumptions \cite{farivar2013branch2 ,huang2016sufficient}, which typically do not hold under constant load scenarios considered here.\vspace{-2.5ex}
\subsection{Inexactness Metrics for SOC Relaxation}
To quantify the inexactness introduced by relaxing the SOC‐surface equalities, we frame \eqref{formulation} as a special case of \emph{Second‐Order Cone Surface Programming (SOCSP)}.
\begin{definition}[Second‐Order Cone Surface Programming (SOCSP)]
A SOCSP is an optimization problem of the form:
    \begin{subequations}
        \begin{align}
        \min\ &\mathbf{c}^\top\mathbf{x}\\
        \textrm{s.t.}\ &A\mathbf{x}=\mathbf{b},\\
         &\left\|\bar{\mathbf{x}}_i\right\|_2=x^0_i, i\in [r]\label{socs}
    \end{align}\label{socsp}
    \end{subequations}\noindent
where \(\mathbf{x}=[\mathbf{x}_1^\top, \mathbf{x}_2^\top, \ldots,\mathbf{x}_r^\top]^\top\), \(\mathbf{x}_i=[x_i^0,\bar{\mathbf{x}}_i^\top]^\top\in \mathbb{R}^{n_i}\) for \(i\in [r]\) are variables and \(\mathbf{c}, A, \mathbf{b}\) are given.
\end{definition}

In this notation, formulation \eqref{formulation} is a SOCSP with variable‐pairs ($S_l$, $[P_l,\;Q_l]^\top$) and ($\frac{W_i + \Phi_l}{2}$, $[S_l,\;\frac{W_i - \Phi_l}{2}]^\top$)
for each branch $l\in\mathcal{L}$. The nonconvexity of SOCSP stems solely from the surface constraints \eqref{socs}.
Relaxing each equality in \eqref{socsp} to an inequality, $\|\bar{\mathbf{x}}_i\|_2 \le x^0_i$, yields a standard SOCP, whose solution may no longer satisfy the original surfaces exactly. We now introduce metrics to measure this deviation and tightness of relaxations.

Consider any relaxation of the SOCSP \eqref{socsp} of the form
\begin{equation}
    \min\;\{\mathbf{c}^\mathsf{T}\mathbf{x} : A\mathbf{x}=\mathbf{b},\;\mathbf{x}\in\mathcal{C}\},
    \label{socr_gen}
\end{equation}
with feasible region 
\(\mathcal{X}_{\rm rel}=\{\mathbf{x}:A\mathbf{x}=\mathbf{b},\;\mathbf{x}\in\mathcal{C}\}\),
where \(\mathcal{C}\) is any superset of the original SOC surfaces.
\begin{definition}[Conic Error and \(\varepsilon\)-Feasibility]
   For each block \(i=1,\dots,r\) and any \(\mathbf{x}\in\mathcal{X}_{\rm rel}\), let
\[
\delta_i(\mathbf{x})
=\|\bar{\mathbf{x}}_i\|_2^2 - (x^0_i)^2
\]
denote the \emph{conic error} for the \(i\)th SOC surface. Then let
\[
\delta_i^{\rm abs}(\mathbf{x})
=\bigl|\delta_i(\mathbf{x})\bigr|,\qquad
\delta_i^{\rm rel}(\mathbf{x})
=\frac{\bigl|\delta_i(\mathbf{x})\bigr|}{(x^0_i)^2},
\]
denote \emph{absolute conic error} and \emph{relative conic error}, respectively. The vectors
\(\boldsymbol{\delta}^{\rm abs}(\mathbf{x})\) and \(\boldsymbol{\delta}^{\rm rel}(\mathbf{x})\)
collect these values over \(i=1,\dots,r\). 
The relaxation \eqref{socr_gen} is \(\varepsilon\)\emph{-feasible} if \emph{every}
\(\mathbf{x}\in\mathcal{X}_{\rm rel}\) satisfies
\[
\bigl\|\boldsymbol{\delta}^{\rm rel}(\mathbf{x})\bigr\|_\infty
=\max_{i=1,\dots,r}\delta_i^{\rm rel}(\mathbf{x})
\;\le\;\varepsilon.
\]
\end{definition}


In practice we rely on the relative conic error (with a small $\eta>0$ in the denominator for numerical stability) as our primary inexactness metric. This measure is geometrically analogous to the normalized distance from a candidate solution to the SOC surface, directly quantifying its primal deviation from the equality constraint. Consequently, imposing a smaller tolerance \(\varepsilon\) in the \(\varepsilon\)-feasibility criterion enforces a tighter relaxation. In our formulation, two variants of the SOC relaxation lead to distinct error expressions: 4-D cone:\(\delta_l^{\mathrm{rel}}=|P_l^2 + Q_l^2 - \Phi_l\,W_i|/(\tfrac{W_i + \Phi_l}{2})^{2}\), and decomposed cones:  \(\delta_{l1}^{\mathrm{rel}}
    = |P_l^2 + Q_l^2 - S_l^2|/S_l^2
    ,
    \delta_{l2}^{\mathrm{rel}}
    = |S_l^2 - \Phi_l\,W_i|/
           (\tfrac{W_i + \Phi_l}{2}\bigr)^{2}\).

\section{Piecewise Relaxations and Efficient Formulations}\label{sec:pr}
This section introduces piecewise relaxations of the 3-D SOC surface and presents a compact mixed integer implementation via the rotation‐and‐fold (R\&F) strategy.\vspace{-2ex}
\subsection{Pyramidal Approximation and R\&F Strategy}
The Pyramidal Approximation (PA), introduced by Zhou~\textit{et al.}~\cite{zhou2020pyramidal}, constructs an inner approximation of the 3D SOC surface~\eqref{socp_new} by inscribing a polyhedral cone whose planar facets yield linear inequalities. To encode this structure efficiently, PA employs a rotation-and-fold (R\&F) scheme, which also serves as the foundation of our proposed framework. We begin by summarizing the core construction of PA.
\begin{figure}\vspace{-3ex}
    \centering
    \includegraphics[width=1.3in]{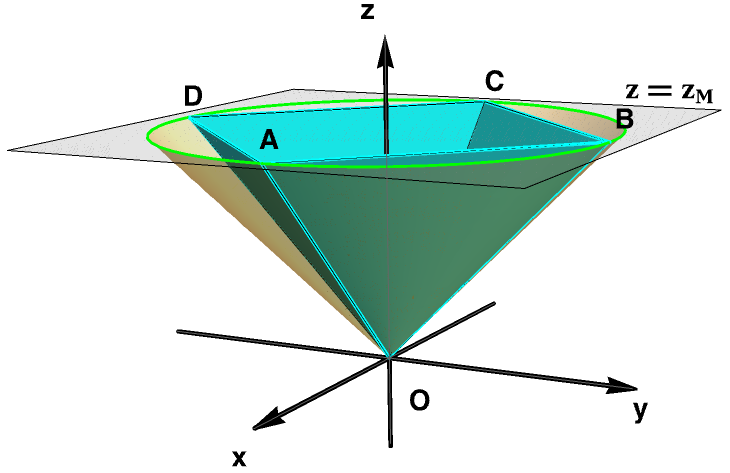}
    \caption{Square-pyramidal
 approximation of the SOC surface \(\sqrt{x^2+y^2}=z\)}
    \label{fig:pa}\vspace{-3ex}
\end{figure}
\subsubsection{Basic idea}
Given the SOC surface \(\sqrt{x^2+y^2}=z\) and an upper bound \(z\le z_M\), PA proceeds as follows:
\begin{itemize}
    \item Select $N$ evenly spaced points $V_1,\dots,V_N$ on the circle 
        $x^2 + y^2 = z_M^2$ in the plane $z = z_M$.
        \item Form the pyramid apex at the origin $O$, and connect $O$ to each pair $(V_i,V_{i+1})$, yielding $N$ triangular facets, where $V_{N+1}=V_1$.
        \item Replace the original SOC surface with these facets.
\end{itemize}
For the case $N=4$, this yields a square-pyramidal
 approximation as Fig. \ref{fig:pa} shows. Every point $P$ in the pyramidal facets admits the convex combination 
\begin{equation}
  \overrightarrow{OP}
  = \sum_{n=1}^N \omega_n\,\overrightarrow{OV_n},
  \label{eq:PA-comb}
\end{equation}
where the coefficients \( \omega_n \) satisfy:
\begin{equation}
  \sum_{n=1}^N \omega_n \le 1,\quad
  \omega_n \ge 0,\quad \omega_n \in \text{SOS2}.
  \label{eq:PA-weights}
\end{equation}
The SOS2 (Special Ordered Set Type 2) condition~\cite{wolsey2020integer} ensures that only two adjacent weights are nonzero, thereby encoding a triangle in the pyramid. By introducing binary variables to model SOS2 constraints, PA yields a MILP formulation whose accuracy improves with increasing \( N \).
\subsubsection{Rotation-and-Fold Strategy}
\begin{figure*}\centering\vspace{-3ex}
    \subfloat[\footnotesize Original $z$-plane]{\includegraphics[scale=.23]{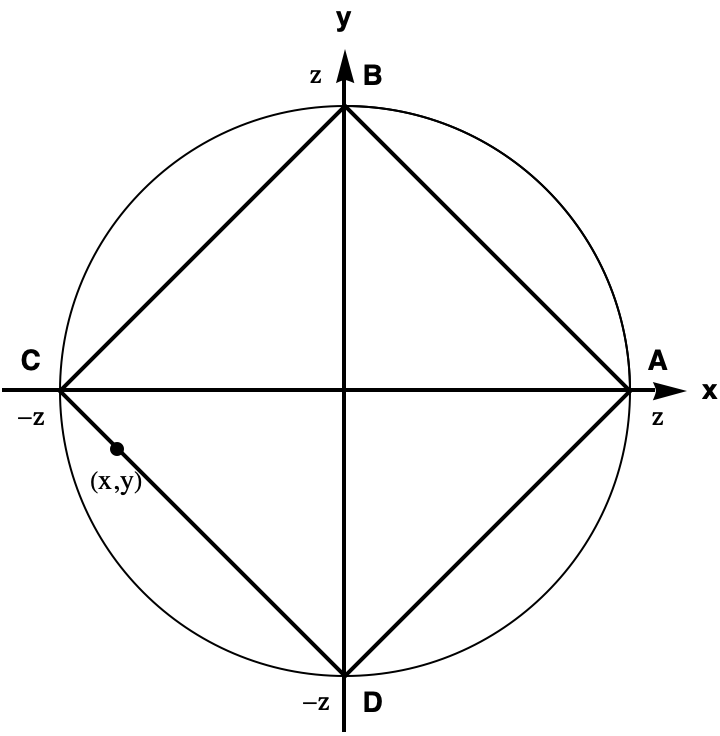}}\quad\quad\quad
  \subfloat[\footnotesize Axial-symmetry]{\includegraphics[scale=.23]{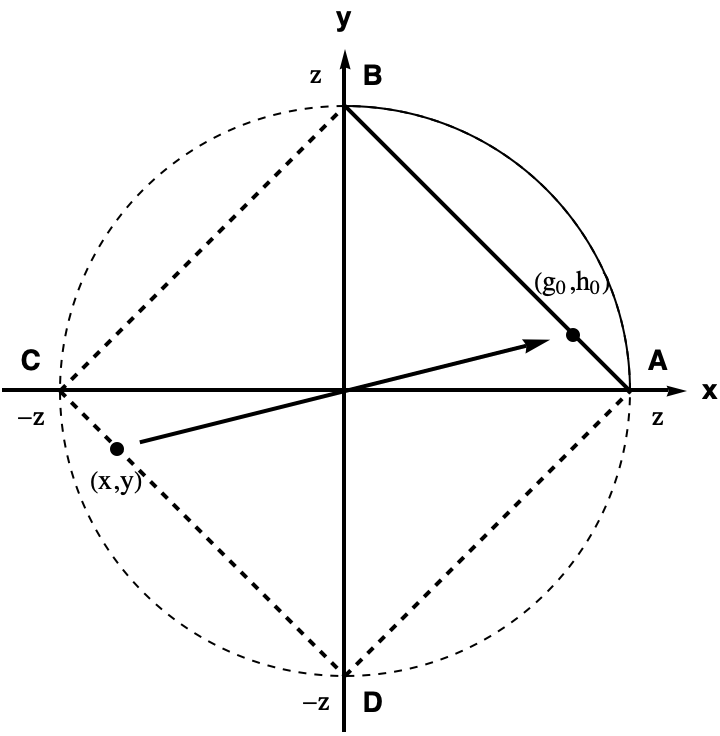}}\quad\quad\quad\quad
  \subfloat[\footnotesize Rotation]{\includegraphics[scale=.18] {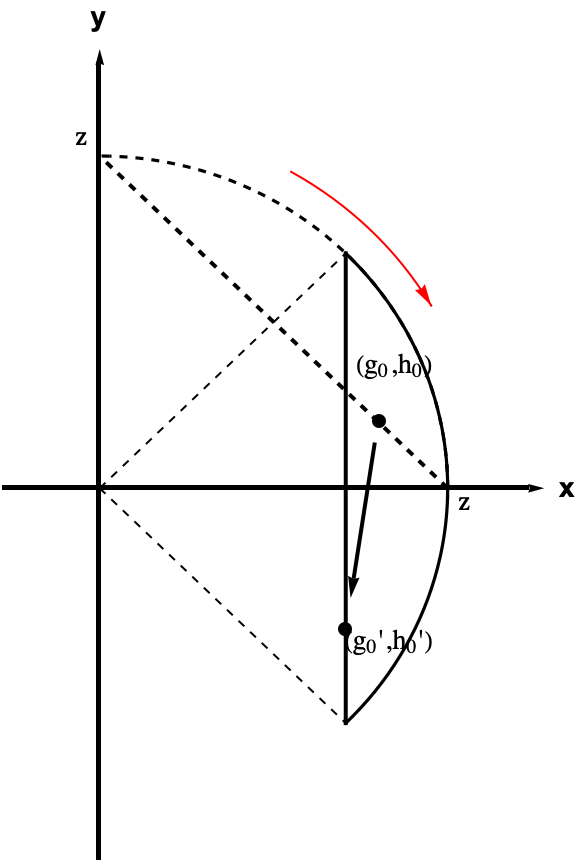}}\quad\quad\quad\quad
  \subfloat[\footnotesize Fold]{\includegraphics[scale=.18] {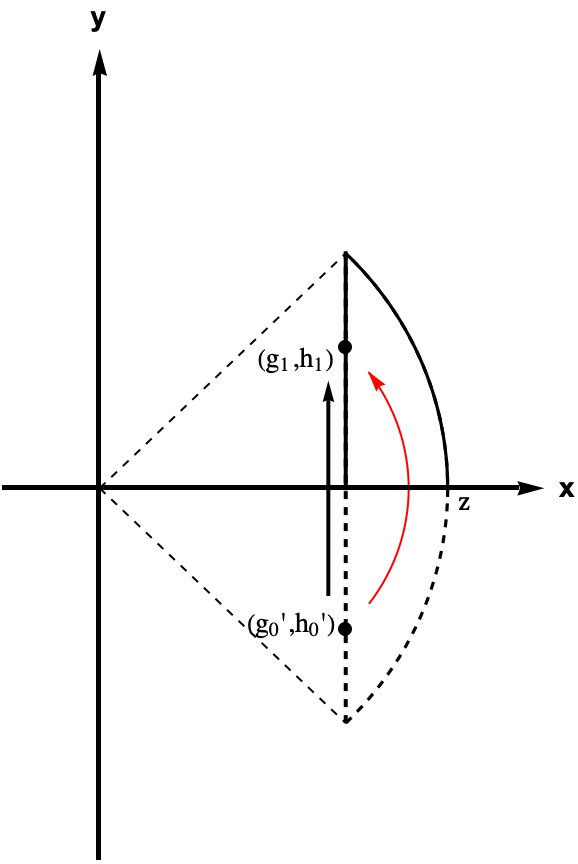}}
  \caption{R\&F strategy in PA}\label{rnf}\vspace{-3ex}
\end{figure*}
The direct MILP encoding of an $N$‐facet pyramidal approximation requires $\mathcal{O}(N)$ binary variables and constraints, which becomes computationally intractable as $N$ grows \cite{zhou2020pyramidal}. To address this, PA adopts a Rotation-and-Fold (R\&F) strategy inspired by Ben-Tal and Nemirovski~\cite{ben2001polyhedral}, which exploits the radial symmetry of the 3D SOC surface to reduce formulation size to \( \mathcal{O}(\log N) \).

We first define the feasible region of the direct pyramid formulation.
\begin{definition}[\(\Pi_K\): Direct Pyramid]\label{def:PiK}
Let $\Pi_K$ be a set of the regular $2^{K+1}$‐sided pyramidal facets inscribed in 
\(\sqrt{x^2+y^2}\le z\le z_M\). Its direct MILP representation is
\begin{equation}
  \Pi_K
  := \Bigl\{(x,y,z)\in\mathbb{R}^2\times[0,z_M]\;:\;
  \eqref{eq:PA-weights},\;\eqref{vceq},\;\bm\omega\in\mathrm{SOS2}\Bigr\},
  \label{eq:PiK-direct}
\end{equation}
where \eqref{eq:PA-weights} are the convex‐combination constraints with \(N=2^{K+1}\), and \eqref{vceq} represent facets in every \(z\)-plane:
\begin{equation}
   \left[\begin{array}{l}
x \\
y
\end{array}\right]=\sum_{n=1}^{2^{K+1}} \omega_n\left[\begin{array}{c}
z \cos \varphi_n \\
z \sin \varphi_n
\end{array}\right]\label{vceq}
\end{equation}
with constants $\varphi_n = \frac{2(n-1) \pi}{2^{K+1}}$. 
\end{definition}

On each fixed \(z\)-plane, the R\&F strategy encoding of the $2^{K+1}$-sided pyramid \(\Pi_K\) proceeds in following three stages:  
\begin{itemize}
    \item {\bf Axial symmetry}\ Reflect all points into the first quadrant:
    \begin{equation}
        \left\{
    \begin{aligned}
        &g_0 = |x|,\\
        &h_0 = |y|,
    \end{aligned}\right.\label{rf1}
    \end{equation}
(Folds quadrants II–IV into quadrant I; see Fig.~\ref{rnf}a–b.)
    \item {\bf Iterated rotation and fold}\ For $k=1,\dots,K$, rotate by angle $\theta_k=\pi/2^{k+1}$ and fold:
    \begin{equation}
        \left\{\begin{aligned}
& g_k=\cos \theta_k g_{k-1}+\sin \theta_k h_{k-1}, \\
& h_k=\left|-\sin \theta_k g_{k-1}+\cos \theta_k h_{k-1}\right|,
\end{aligned}\right.\quad k\in[K]. \label{rfmap}
    \end{equation}
    Each rotation by $\theta_k$ and subsequent fold (absolute‐value) operation halves the remaining angular sector and facet (Fig.~\ref{rnf}c–d).
    \item {\bf Final linear facet} After $K$ iterations, the residual sector has angle $\theta_K$, and only a line segment of the facet remains which can be formulated by simple linear constraints:
\begin{equation}
  g_K = z\,\cos\theta_K,
  \qquad
  0 \le h_K \le z\,\sin\theta_K.
  \label{rffinal}
\end{equation}
\end{itemize}

The absolute-value operation in \eqref{rf1} and \eqref{rfmap} are linearized via the big-$M$ method, which replaces the constraint $Y=|X|$ with:
\begin{equation}\label{big-m}
\left\{\begin{aligned}
& X=M\left(\omega_1-\omega_2\right), \\
& Y=M\left(\omega_1+\omega_2\right), \\
& 0 \leq \omega_1 \leq \beta, \\
& 0 \leq \omega_2 \leq 1-\beta,
\end{aligned}\right.
\end{equation}
where $\omega_1$ and $\omega_2$ are continuous variables, $\beta$ is a binary variable indicating the sign of $X$. The constant $M$ is chosen as a valid upper bound on $|X|$ to ensure numerical stability and tighten the formulation. With this reformulation, the R\&F encoding becomes a pure MILP approach.
\begin{theorem}[Equivalence of R\&F and Direct Pyramid]\label{thm:rnf}
For any integer \(K\ge1\) and \((x,y,z)\in\mathbb{R}^2\times[0,z_M]\), $ (x,y,z)\in\Pi_K $ if and only if \eqref{rf1}, \eqref{rfmap}, and \eqref{rffinal} all hold, with each $|\cdot|$ replaced by the big-\(M\) system \eqref{big-m}.
\end{theorem}

Theorem \ref{thm:rnf} shows that R\&F strategy yields an equivalent mixed integer encoding of the $2^{K+1}$-facet pyramid using only $\mathcal O(K)$ binaries and linear constraints, compared to the $\mathcal O(2^{K+1})$ required by direct SOS2 formulations, achieving identical approximation tightness with a logarithmic formulation size.\vspace{-2ex}
\subsection{Piecewise Relaxations of SOC Surface}
Although PA provides a tighter approximation than standard SOCP, its approximation nature can render the formulation infeasible. 
To see this, let $\mathcal{C}$ be a certain 3-D SOC surface defined by \eqref{cone1} for branch \(l\) and $\mathcal{P}$ be the corresponding polyhedral projection in \((P_l,Q_l,S_l)\) space defined by all linear constraints in \eqref{formulation}. Any feasible OPF must satisfy \(\mathcal{P}\cap\mathcal{C}\neq\varnothing\). However, PA replaces \(\mathcal{C}\) by its inscribed pyramid \(\Pi_K\). It is possible that \(\Pi_K\cap \mathcal{C}=\varnothing\) causing the PA-based formulation to become infeasible. 

Unlike approximations, relaxations preserve feasibility provide valid dual bounds for the original problem. To this end, we propose two piecewise relaxations of the 3D SOC surface: \emph{Pyramidal Relaxation} (PR) and \emph{Quasi-Pyramidal Relaxation} (QPR). 

\begin{figure}\centering\vspace{-3ex}
    \subfloat[]{\includegraphics[scale=.2]{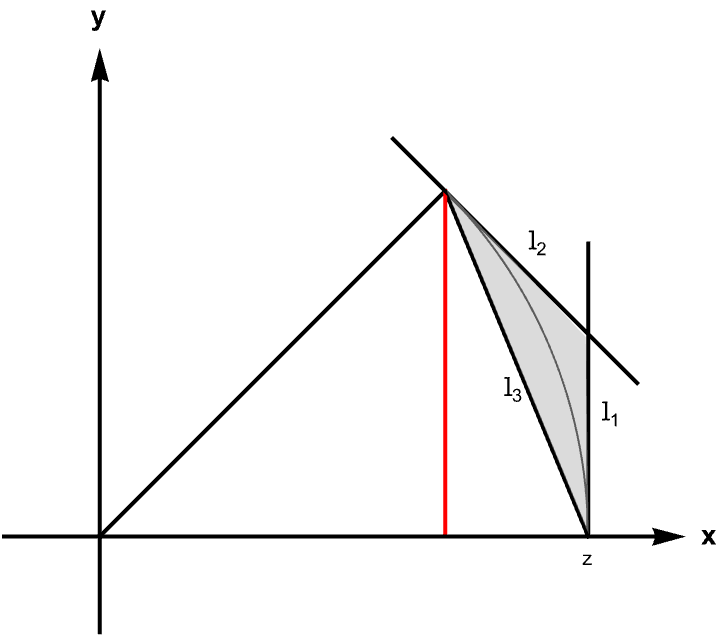}}\quad\quad\quad
  \subfloat[]{\includegraphics[scale=.2]{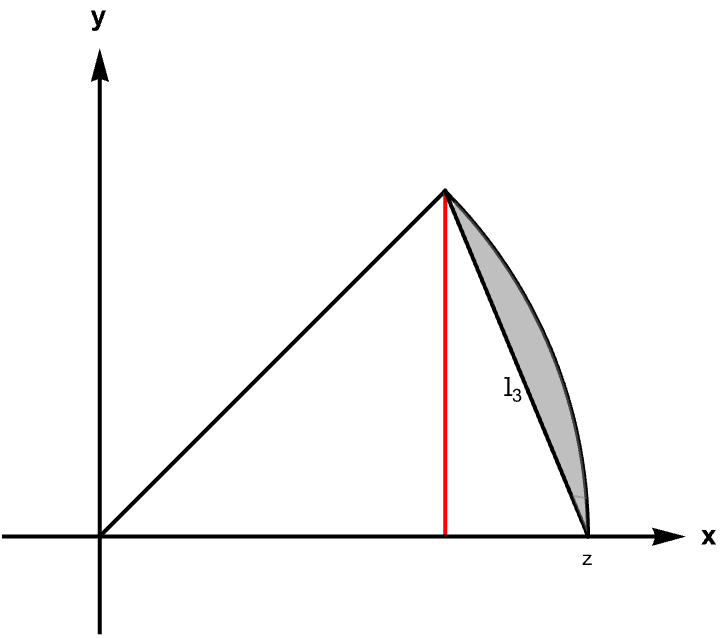}}
  \caption{Final relaxation regions in the $z$‐plane for PR (left) and QPR (right).}\label{relax}\vspace{-3ex}
\end{figure}

Both PR and QPR reuse the R\&F steps for efficiency. The axial symmetry \eqref{rf1} and iterated rotation‐and‐fold \eqref{rfmap} stages remain unchanged.  The difference lies in the final step:

\paragraph{Pyramidal Relaxation (PR)} PR replaces the single linear facet \eqref{rffinal} by a triangular region in the $z$‐plane bounded by lines $l_1,l_2,l_3$ as illustrated in Fig.~\ref{relax}a: 
\begin{equation}
\left\{\begin{aligned}
& g_K \leq z, \\
& g_K \cos \theta_K+h_K \sin \theta_K \leq z,\\
& z \cos \theta_{K+1} \leq g_K \cos \theta_{K+1}+h_K \sin \theta_{K+1},
\end{aligned}\right.\label{eqpr}
\end{equation}
The first two inequalities in \eqref{eqpr} are \emph{outer cuts}, lying outside the cone, while the third is an \emph{inner cut}. Geometrically, PR defines a union of \(N = 2^{K+1}\) tetrahedra that include the SOC surface in 3D.
\paragraph{Quasi‐Pyramidal Relaxation (QPR)}  
QPR modifies the PR formulation by retaining only the inner cut:
\begin{equation}
    z \cos \theta_{K+1} \leq g_K \cos \theta_{K+1}+h_K \sin \theta_{K+1},\label{eq:QPR}
\end{equation}
and preserving the original SOC constraint:
\begin{equation}
    \sqrt{x^2+y^2}\le z. \label{eq:QPR-socp}
\end{equation} 
This produces a circular-segment region on each \(z\)-plane (Fig.~\ref{relax}(b)), corresponding in 3D to a union of circular-segmental cone, or ``quasi-pyramids".

For notational clarity, we collect the three $K$‐stage R\&F variants below. Each keeps the base objective $\mathbf{c}^\top\mathbf{x}$ and linear constraints and replaces each SOC block constraint by the indicated R\&F equations:

\begin{definition}[$K$-Stage R\&F Approximations and Relaxations]\label{def:pa,pr,qpr}
For any integer \(K\ge0\), consider the 3D SOCSP \eqref{socsp} with each block \(\mathbf{x}_i\in\mathbb{R}^3\). Denote by
$\mathcal{X}_{\text{LP}}
=\{\mathbf{x} : A\mathbf{x}=\mathbf{b}\}$
the set of the linear constraints. The feasible regions of the three R\&F-based formulations are:
\begin{itemize}
  \item {\bf $K$‐PA (Pyramidal Approximation)} 
  \[
    \mathcal{X}_{K\text{-PA}}
    =\mathcal{X}_{\text{LP}}\cap \bigl\{\mathbf{x}:
      \eqref{rf1},\;\eqref{rfmap},\;\eqref{rffinal} \text{ hold }\forall i\in[r] \bigr\}.
  \]
  \item {\bf $K$‐PR (Pyramidal Relaxation)} 
  \[
    \mathcal{X}_{K\text{-PR}}
    =\mathcal{X}_{\text{LP}}\cap\bigl\{\mathbf{x}:
      \eqref{rf1},\;\eqref{rfmap},\;\eqref{eqpr}\text{ hold }\forall i\in[r] \bigr\}.
  \]
  \item {\bf $K$‐QPR (Quasi-Pyramidal Relaxation)} 
  \[
    \mathcal{X}_{K\text{-QPR}}
    =\mathcal{X}_{\text{LP}}\cap\bigl\{\mathbf{x}:
      \eqref{rf1},\;\eqref{rfmap},\;\eqref{eq:QPR},\;\eqref{eq:QPR-socp}\text{ hold }\forall i\in[r]
      \bigr\}.
  \]
\end{itemize}
Here, \(K\) is the number of rotation‐and‐fold iterations in \eqref{rfmap}. All absolute‐value terms are linearized via the big-\(M\) method \eqref{big-m}.  
\end{definition}

Note that \(K \ge 1\) is required for \(K\)-PA to define the final facet via~\eqref{rffinal}, whereas \(K \ge 0\) suffices for PR and QPR. Notably, the outer facet of \((K+1)\)-PA coincides with the inner boundary of \(K\)-PR and \(K\)-QPR, leading to the set inclusion: \[\mathcal{X}_{(K+1)\text{-PA}}\subseteq\mathcal{X}_{K\text{-QPR}}\subseteq\mathcal{X}_{K\text{-PR}}.\] Thanks to the R\&F strategy, all three formulations require only \( \mathcal{O}(K) \) binary and linear constraints to partition each cone surface into \(2^{K+2}\) parts.
From a modeling perspective, PA and PR are pure MILP formulations, whereas QPR is a Mixed Integer Second‐Order Cone Programming (MISOCP) due to its retained SOC constraint \eqref{eq:QPR-socp}. Although the nonlinear cone may slow down the solver, QPR provides the tightest relaxation, as quantified below.
\begin{theorem}\label{thm:eps}
    For any integer \(K\ge0\) and any 3D SOCSP,
    \begin{itemize}
        \item \(K\)-PA is \(\sin^2\theta_K\)-feasible,
      \item \(K\)-PR is \(\tan^2\theta_{K+1}\)-feasible,
      \item \(K\)-QPR is \(\sin^2\theta_{K+1}\)-feasible,
    \end{itemize}
    where \(\theta_K = \pi/2^{K+1}\).
\end{theorem}
In practice, one selects \(K\) to meet a target error tolerance. For example, to enforce relative conic error below 1\%, one may use 4-PA, 3-PR or 3-QPR.  Moreover, different branches can be assigned different \(K\) values to reflect their varying sensitivity to SOC relaxation. As \(K\) increases, all three formulations converge asymptotically to the feasible region and thus the global optimum of the original SOCSP.

\vspace{-1.5ex}
\section{Dynamic Relaxation Algorithm Framework}\label{sec:drf}
In this section, we propose a dynamic relaxation algorithm framework that integrates the R\&F‐based PR and QPR into a branch‐and‐cut algorithm, making it directly compatible with commercial solvers. By leveraging the nested structure of these relaxations, the framework dynamically generates violated cuts during the solve, progressively tightening the SOC relaxation without model reconstruction. To further enhance efficiency and solution quality, we incorporate warm-start initialization and post-processing heuristics.\vspace{-3ex}
\subsection{Motivation}
\begin{figure}\centering\vspace{-3ex}
    \subfloat[PA]{\includegraphics[scale=.3]{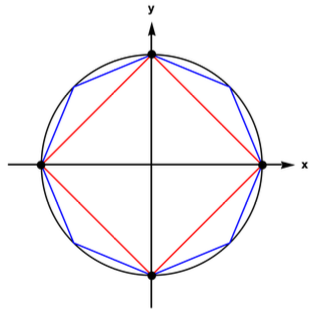}}\quad \quad\quad
  \subfloat[PR]{\includegraphics[scale=.3]{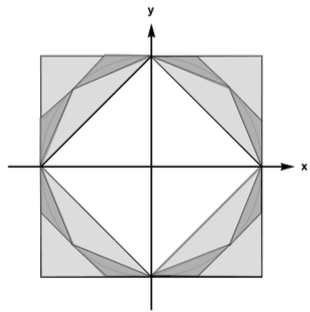}}\quad \quad\quad
  \subfloat[QPR]{\includegraphics[scale=.3]{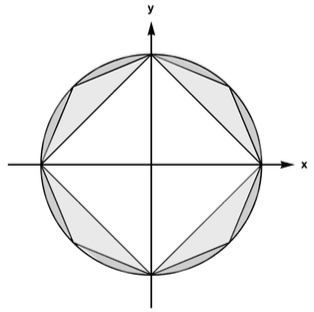}}
  \caption{Static approximation vs.\ dynamic relaxations on \(z\)-plane.
(a) \textbf{PA}: the red (1-PA) and blue (2-PA) polygons are non-nested. 
(b) \textbf{PR}: 0-PR (entire gray) contains 1-PR (dark gray).
(c) \textbf{QPR}: 0-QPR (entire gray) contains 1-QPR (dark gray).}
  \label{mot}\vspace{-3ex}
\end{figure}

Although PA can produce tight approximations, it is inherently static: its feasible regions for different depths \(K\) do not nest. For instance, simply adding the R\&F mapping and the facets of \(\mathcal{X}_{2\text{-PA}}\) to \(\mathcal{X}_{1\text{-PA}}\) yields disconnected elements in the \(z\)-plane (Fig.~\ref{mot}a), offering little benefit. In contrast, both PR and QPR form nested relaxations: 
\[
\mathcal{X}_{K\text{-PR}}\;\supseteq\;\mathcal{X}_{(K+1)\text{-PR}}
\quad\text{and}\quad
\mathcal{X}_{K\text{-QPR}}\;\supseteq\;\mathcal{X}_{(K+1)\text{-QPR}},
\]
for any \(K\ge0\), allowing progressive refinement by simply appending additional cuts from deeper R\&F iterations (Fig.~\ref{mot}b–c). This nesting property enables efficient integration with branch-and-cut solvers.

A second key insight concerns constraint activeness. Although each \(K\)-stage PR or QPR nominally introduces \(\mathcal{O}(K)\) cuts per branch, only a small subset of them binds at optimality. This motivates a lazy cut-generation strategy: we initialize the model with a coarse relaxation, and iteratively add only those cuts violated by incumbent solutions. This reduces problem size, improves solver speed, and ensures tightness only where necessary.\vspace{-2.5ex}
\subsection{Algorithm Outline}
We now describe the proposed \emph{dynamic relaxation framework}, which embeds PR or QPR within a branch-and-cut algorithm via specialized cut-generation mechanism. The framework begins with a coarse \(K_{\text{init}}\)-stage formulation (either PR or QPR) and, optionally, a warm-start incumbent. During the solve, each integer-feasible solution is checked for \(\varepsilon\)-feasibility. If the relative conic error exceeds the tolerance, new R\&F cuts from deeper stages are added on the fly, tightening the relaxation locally without reconstructing the model. This solve–verify–cut loop continues until standard solver termination criteria are met (e.g., optimality gap or time limit).
To improve solution quality, a post-processing heuristic may be applied after convergence. The overall framework is illustrated in Fig.~\ref{fig:flowchart},
and full pseudocode is provided in Appendix \ref{appendix-alg}.

We refer to the PR-based instantiation as the \emph{Dynamic Pyramidal Relaxation} (DPR) method, and the QPR-based variant as the \emph{Dynamic Quasi-Pyramidal Relaxation} (DQPR). When the cut depth is capped at level \(K\), we denote the methods as \(K\)-DPR and \(K\)-DQPR, respectively. The following proposition establishes consistency with their static counterparts.
\begin{proposition}\label{prop:consist}
For any integer \(K \ge 0\), the \(K\)-DPR (resp.\ \(K\)-DQPR) algorithm produces an optimal solution to the static \(K\)-PR (resp.\ \(K\)-QPR) formulation.
\end{proposition}

\begin{figure}\vspace{-3ex}
    \centering
    \includegraphics[width=3in]{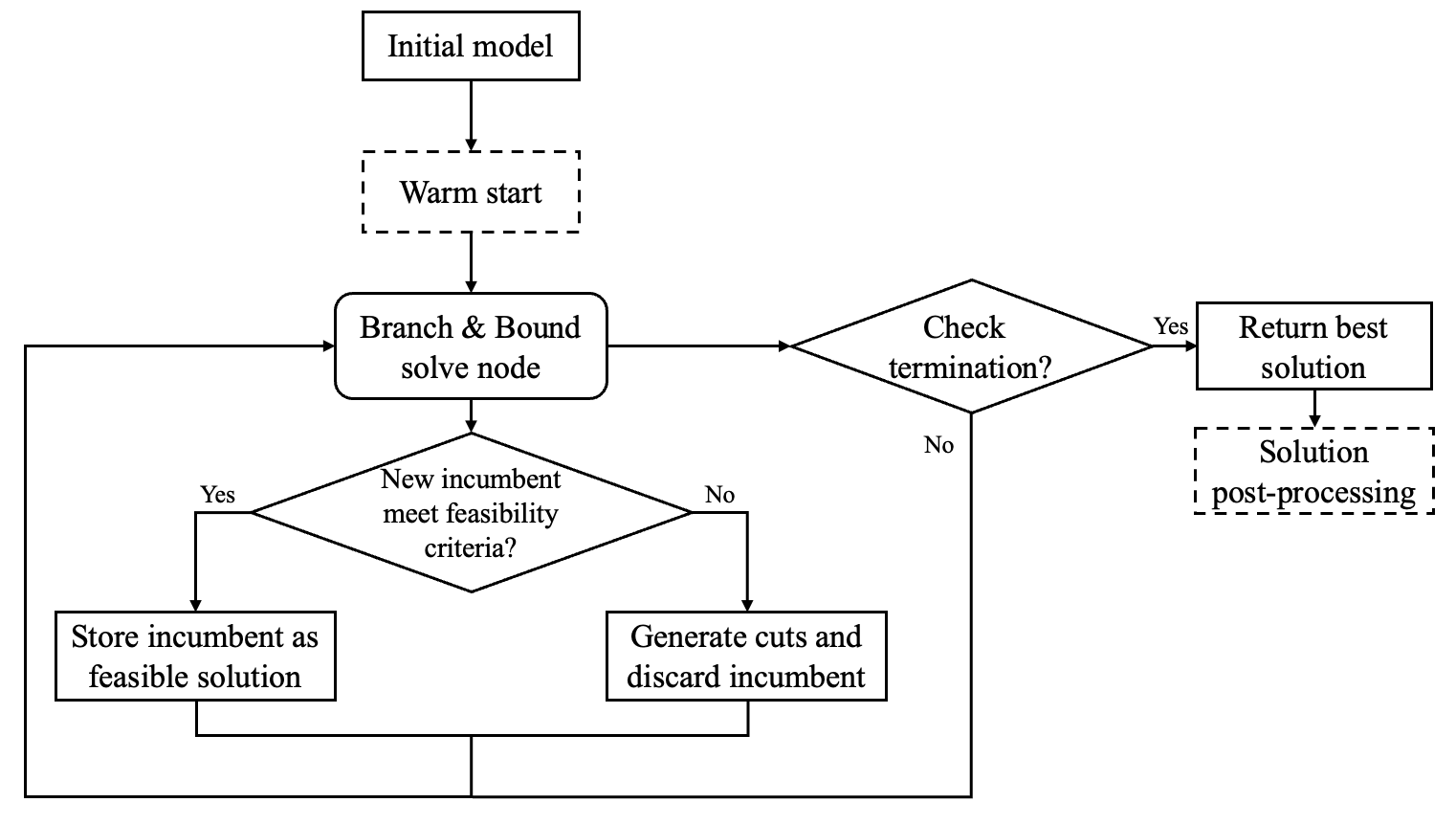}
    \caption{Overview of the dynamic relaxation framework.}
    \label{fig:flowchart}\vspace{-3ex}
\end{figure}\vspace{-2.5ex}
\subsection{Algorithm Details} \label{subsec:ald_det}
\subsubsection{Cut generation}
Cut generation is the core mechanism of the dynamic framework. At each branch-and-cut node, candidate solutions are tested for \(\varepsilon\)-feasibility. If a violation is detected on a branch \(l\), targeted cuts are generated to tighten the corresponding relaxation. 
In PR, two types of R\&F-based cuts are used: \emph{outer cuts}, which intersect the current feasible region via additional linear inequalities, and \emph{inner cuts}, which partition the cone into wedges and require binary variables to encode the resulting union. QPR uses only inner cuts, as the original SOC constraint already provides an outer envelope.
Due to this structural distinction, the framework applies cut generation selectively: outer cuts are added incrementally without introducing new binaries, while inner cuts are added more conservatively and only at the minimum depth required to eliminate the violated solution. We now describe each procedure in detail.

\subsubsection*{Inner-cut generation}
If a candidate solution \(\mathbf{x}^*\) violates the \(\varepsilon\)-feasibility condition while remaining inside the cone on branch \(l\), we trigger the inner-cut routine. Let \(k_0\) denote the current R\&F depth for that branch. We then search for the smallest \(k_1 > k_0\) such that the corresponding inner-cut inequality,
\begin{equation}\label{incut}
    z\cos\theta_{k_1+1}
\;\le\;
g_{k_1}\cos\theta_{k_1+1} + h_{k_1}\sin\theta_{k_1+1}
\end{equation}
is violated by \(\mathbf{x}^*_l\). To exclude this solution, we append the R\&F mappings for levels \(k_0+1\) through \(k_1\), along with the \(k_1\)-level inner cut. This introduces binary variables, but only at the minimal depth needed.
Figure~\ref{op}a illustrates this update in the \(z\)-plane: the previous feasible wedge is narrowed by new inner cuts, eliminating \(\mathbf{x}^*_l\) from the relaxed region.

\subsubsection*{Outer-cut generation}
\begin{figure}\centering\vspace{-3ex}
    \subfloat[]{\includegraphics[scale=.2]{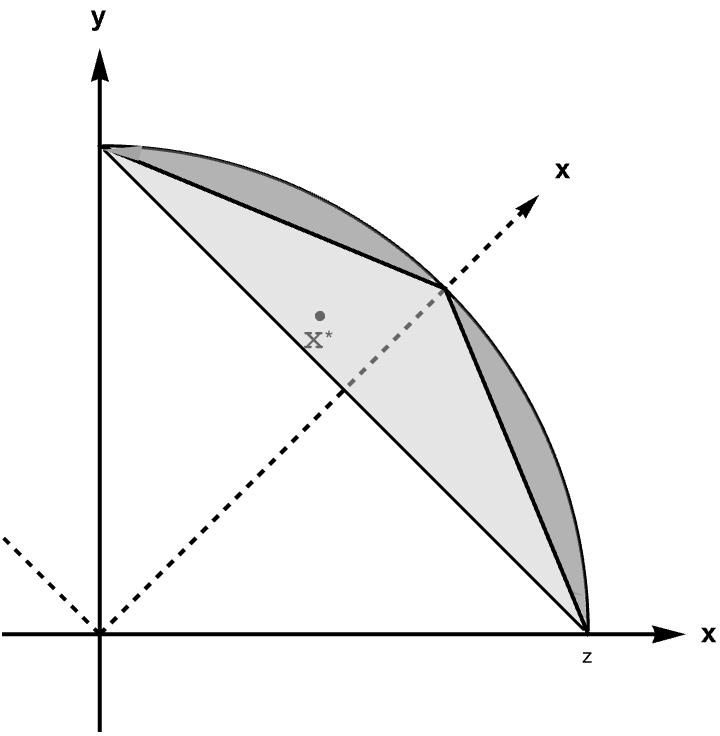}}\quad\quad\quad
  \subfloat[]{\includegraphics[scale=.2]{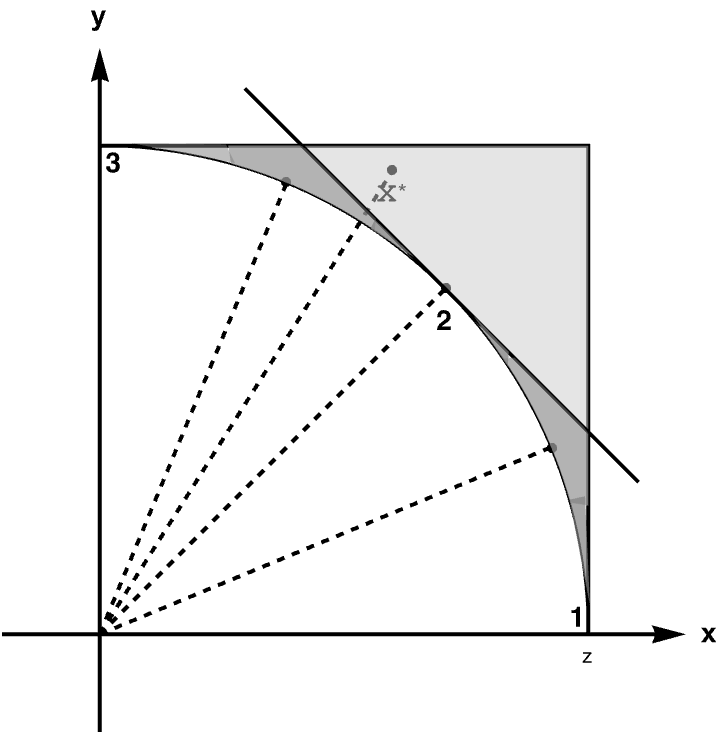}}
  \caption{Illustration of inner-cut (a) and outer-cut (b) generation in the \(z\)–plane. The entire shading area is the
previous region. (a) Dashed axes mark the new R\&F rotation; the dark shading
region is the refined region after adding the $k_1$–inner cut that excludes the
violated solution. (b) After adding new tangent cuts at adjacent angles, the
dark-shaded region excludes the infeasible solution.}\label{op}\vspace{-3ex}
\end{figure}

Outer cuts are used only in DPR to exclude infeasible points outside the cone. When a candidate solution \(\mathbf{x}^*\) lies outside the cone and fails \(\varepsilon\)-feasibility on branch \(l\), we identify the first outer cut level \(k > k_0\) where the violation occurs:
\[
\delta_{k+1}^{\mathrm{out}} < \delta_l(\mathbf{x}^*_l) \le \delta_k^{\mathrm{out}}.
\]
Let \(\theta^* = \arg(\mathbf{x}_l^*)\). We then determine the two adjacent angles \(\psi_1, \psi_2\) at level \(k\) that bracket \(\theta^*\), and add the corresponding tangents:
\begin{equation}
    \left\{
\begin{aligned}
& g_{k_0} \cos \psi_1 + h_{k_0} \sin \psi_1 \le z, \\
& g_{k_0} \cos \psi_2 + h_{k_0} \sin \psi_2 \le z.
\end{aligned}
\right.\label{otc}
\end{equation}
These linear constraints cut off \(\mathbf{x}_l^*\) without introducing new binaries. Fig.~\ref{op}b shows the effect: the new tangents (e.g., at point 2) shrink the outer region, excluding the violated point while preserving formulation simplicity.

\subsubsection*{R\&F mapping refinement}
To streamline outer-cut generation, we enforce a single linear constraint \(g_k \le z\) at each R\&F depth. The following proposition shows that these inequalities are sufficient to imply the full outer-facet inequality at the final depth.
\begin{proposition}\label{thm:bound}
Let \(g_0, h_0 \in [0, z]\), and define the R\&F updates for \(k = 1, \dots, K\):
 \begin{equation}\label{thmeq:rnfmap}
      \left\{
  \begin{aligned}
    g_k &= \cos\theta_k\,g_{k-1} + \sin\theta_k\,h_{k-1},\\
    h_k &= \bigl|-\sin\theta_k\,g_{k-1} + \cos\theta_k\,h_{k-1}\bigr|,\\
    g_k &\le z,
  \end{aligned}\right.
  \end{equation}
Then the final outer facet inequality
 \begin{equation}
      g_K\cos\theta_K + h_K\sin\theta_K \;\le\; z. \label{thmeq:finalbound}
  \end{equation}
holds automatically.
\end{proposition}

This result allows us to propagate tight outer bounds with minimal overhead: instead of explicitly adding both outer-facet cuts at each level, we need only enforce \(g_k \le z\) during each R\&F step. Before adding any outer cut, we also check for redundancy and omit cuts already implied by previous constraints.

\subsubsection{Warm start}
Providing a high-quality initial solution can significantly accelerate convergence in mixed-integer optimization, especially when early incumbents may be discarded by dynamic cuts in our case. To this end, we solve the AC OPF using the IV formulation~\cite{o2012iv} via the IPM, which yields a suboptimal but feasible AC solution \((\tilde{p}_g, \tilde{q}_g, \tilde{\mathbf{v}}_i, \tilde{\mathbf{i}}_{ij})\), where \(\tilde{\mathbf{v}}_i, \tilde{\mathbf{i}}_{ij}\) are complex voltages and branch power flow, respectively. This solution is mapped to our model variables as follows:
\[
\begin{aligned}
  p_g &= \tilde p_g, &
  q_g &= \tilde q_g, &&
  \forall g\in\mathcal G,\\
  P_{ij} &= \Re\bigl(\tilde{\mathbf v}_i\,\tilde{\mathbf i}_{ij}^*\bigr), &
  Q_{ij} &= \Im\bigl(\tilde{\mathbf v}_i\,\tilde{\mathbf i}_{ij}^*\bigr), &&
  \forall i\to j\in\mathcal L,\\
  S_l &= \sqrt{P_l^2 + Q_l^2}, &
  \Phi_l &= \lvert\tilde{\mathbf i}_l\rvert^2, &&
  \forall l\in\mathcal L,\\
  W_i &= \lvert\tilde{\mathbf v}_i\rvert^2, &&
  & & 
  \forall i\in\mathcal N.
\end{aligned}
\]
Here, \(\Re\{\cdot\}\) and \(\Im\{\cdot\}\) denote real and imaginary parts, and \((\cdot)^*\) denotes complex conjugation. The resulting point satisfies all linear and SOC-surface constraints in our model, ensuring \(\varepsilon\)-feasibility and providing a strong primal bound for the solver.
\subsubsection{Solution post-processing}
\begin{figure}\vspace{-3ex}
    \centering
    \includegraphics[width=1in]{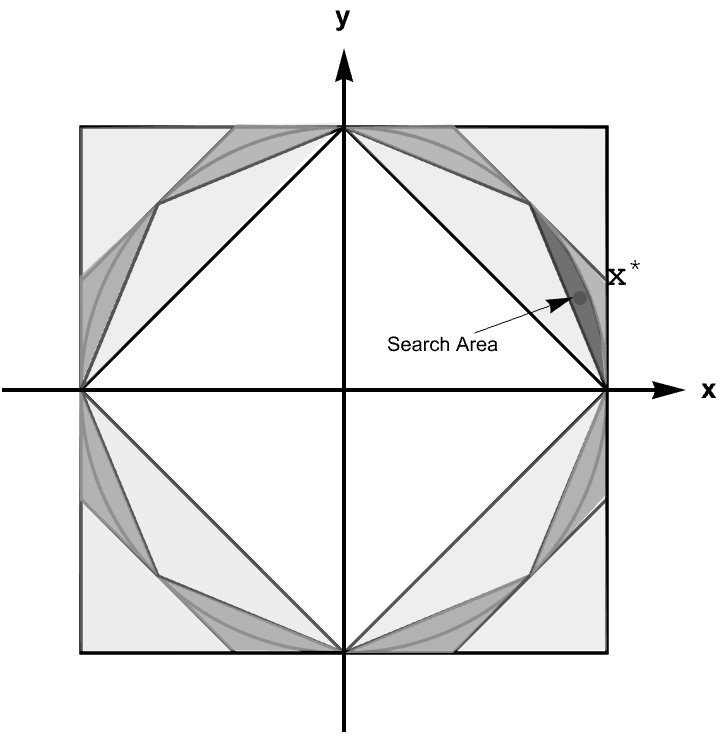}
    \caption{Post‐processing via fixed R\&F region: The shaded union of wedges collapses to a single quasi-pyramid (black) after fixing the R\&F binary pattern and reimposing the SOC constraint.}
    \label{fig:prpp}\vspace{-3ex}
\end{figure}

To further improve solution quality, we apply a simple large-neighborhood search (LNS) as a post-processing step. When DPR terminates, the best incumbent is associated with a fixed pattern of R\&F binaries \(\bm \beta_0^*, \beta_1^*, \ldots, \beta_K^*\), identifying a unique wedge in the relaxation. By fixing these binary variables, the original feasible region—defined as a union of wedge-shaped regions—collapses to a single wedge. Within this wedge, we reimpose the original SOC constraint, converting the problem into a pure SOCP. This subproblem can be solved efficiently and often yields a refined solution with lower relative conic error. Although this post-processed solution is not guaranteed to be optimal for the full MISOCP, it typically improves upon the DPR result at negligible additional cost. As illustrated in Fig.~\ref{fig:prpp}, the final solution lies within a single quasi-pyramid (black), providing a high-quality approximation of the QPR. This framework also allows integration with other heuristics, including machine learning–guided refinement or domain-specific adjustments, highlighting the extensibility of the proposed framework.

\vspace{-2ex}
\section{Numerical Results}\label{sec:nr}
We test the static relaxations (SOCP, PA, PR, QPR) and dynamic methods (DPR, DQPR) on eight instances from PGLib–OPF \cite{babaeinejadsarookolaee2019power}: IEEE 5, 30, 118, 162 DTC, ACTIVSg 200, IEEE 300, ACTIVSg 500, and GOC 793. All experiments were run on a CentOS 7 server with an Intel Xeon Gold 6230R CPU (2.1 GHz) and 128 GB of RAM. Our implementations use Julia 1.9.0 with PowerModels.jl v0.30.0 \cite{coffrin2018powermodels} and JuMP v1.23.3 \cite{Lubin2023}. We solve all MILP/MISOCP formulations using Gurobi 11.0.3 \cite{gurobi} with up to 32 threads, a 3600 seconds time limit, and a 0.1 \% optimality gap; all other solver options are left at their defaults.  For warm starts, we solved the IV‐formulation using IPOPT \cite{wachter2006implementation} and pass its feasible solution to Gurobi. Throughout our experiments, all dynamic runs begin with \(K_{\mathrm{init}}=0\), and cuts are generated via Gurobi callbacks.

The formulations and methods we compare are as follows.
\begin{itemize}
    \item  \textbf{SOCP:}  The classical 4-D SOC relaxation \eqref{socp_proto} from \cite{farivar2013branch}
    \item \textbf{PA}: Pyramidal approximation from \cite{zhou2020pyramidal}
    \item \textbf{PR}: Pyramidal relaxation formulation (static)
    \item \textbf{QPR}: Quasi-pyramidal relaxation formulation (static)
    \item \textbf{DPR}: Dynamic pyramidal relaxation method
    \item \textbf{DQPR}: Dynamic quasi-pyramidal relaxation method
\end{itemize}
To unify comparison, we also parameterize all methods by the total number of piecewise regions \(N\). For PA, \(N = 2^{K+1}\); for PR and QPR, \(N = 2^{K+2}\). 
\vspace{-3ex}
\subsection{Static Formulation Comparison}\label{res1}
We begin by comparing the performance of the static relaxations on a representative test case.
Table \ref{tab:PA‐feas} shows, for each test instance, the smallest pyramid‐size $N=2^{K+1}$ at which PA becomes feasible.  Whenever $N < N_{\mathrm{crit}}$, the inscribed pyramid $\Pi_K$ fails to intersect the polyhedron, causing PA to declare infeasibility. In practice, this means that using a fixed PA partition can lead to infeasible models unless \(N\) is set large enough in advance.  

Table~\ref{tab:static-118} summarizes results for all static methods on the IEEE 118 case, under coarse and fine partitions (\(N=8\) and \(N=128\)). For each method we record the objective, the maximum relative conic error $\|\bm\delta^{\mathrm{rel}}\|_\infty$ (in \%), the maximum absolute conic error $\|\bm\delta^{\mathrm{abs}}\|_\infty$, the sum of absolute conic error $\|\bm\delta^{\mathrm{abs}}\|_1$, and the total solve time. Note that all \(\bm \delta\) in this section are 4-D conic error.

\begin{table}\vspace{-3ex}
  \centering\scriptsize
  \caption{Critical Pyramid Sizes \(N\) for Which PA is First Feasible}
  \label{tab:PA‐feas}
  \begin{tabular}{l c c c c c c c c}
    \toprule
    Instance & 5 & 30 & 118 & 162 & 200 & 300 & 500 & 793 \\
    \midrule
    PA\,\(N_{\mathrm{crit}}\) & 8 & 8 & 8  & 32  & 8  & 128 & 32 & 16 \\
    \bottomrule
  \end{tabular}\vspace{-2ex}
\end{table}
\begin{table}
\scriptsize
 \setlength{\tabcolsep}{5pt}
  \centering
  \caption{Static‐relaxation Results on IEEE 118 for \(N=8\) and\ \(N=128\) }
  \label{tab:static-118}
  \begin{tabular}{%
      l                
      c                
      ccc                 
      c                
    }
    \toprule
    \multirow{2}{*}{Method} 
      & \multirow{2}{*}{Obj} 
        & \multicolumn{3}{c}{Conic error} 
              & \multirow{2}{*}{Time (s)} \\
    &  & {$\|\bm \delta^{\mathrm{rel}}\|_\infty$}(\%) & {$\|\bm \delta^{\mathrm{abs}}\|_\infty$}&{$\|\bm \delta^{\mathrm{abs}}\|_1$}
         &  \\
    \midrule
    SOCP 
      &  96\,334
      &  8.2  &  57.9   & 103.6
      &  0.2 \\ 
    \midrule
    \multicolumn{6}{c}{\textbf{N=8}} \\[0.5ex]
    PA   
      & 101\,665
      & 23.5   & 68.6 & 191.9
      &  {\bf 1.9} \\ 
    PR   
      &  95\,783
      & 29.5   & 62.8 & 178.7
      &  4.6 \\ 
    QPR  
      &  {\bf 96\,334}
      &  {\bf 8.3}   &  {\bf 60.5}   & {\bf 106.2}
      & 15.3 \\ 
    \midrule
    \multicolumn{6}{c}{\textbf{N=128}} \\[0.5ex]
    PA   
      &  96\,378
      &  0.1   &   0.3   & 0.6
      & 214.8 \\ 
    PR   
      &  96\,425
      &  0.1   &   0.1  & 0.2
      & {\bf 171.5} \\ 
    QPR  
      &  {\bf 96\,412}
      &   \(<\mathbf{0.1}\) &  \(<\mathbf{0.1}\) & \(<\mathbf{0.1}\)
      & 573.2 \\ 
    \bottomrule
  \end{tabular}\vspace{-3ex}
\end{table}
When \(N=8\), PA and PR both exhibit large conic errors and moderate runtimes, whereas QPR—though matching SOCP’s 8 \% relative error—requires about three times longer than PR. In other words, at this coarse partition size only QPR achieves SOCP‐level accuracy, while PA and PR suffer much larger violations. At the finer partition \(N=128\), all three pyramidal models drive the maximum relative violation below 0.1 \% and the maximum absolute violation under 1.  Here PR solves approximately 20 \% faster than PA, and QPR effectively eliminates all errors (all metrics under 0.1 \%), but at a drastically higher runtime. In summary, PR consistently matches or surpasses PA in tightness while using fewer resources, and QPR alone guarantees near‐zero conic error even for moderate \(N\), provided one is willing to pay the computational cost. Consequently, PR offers the best balance of accuracy and speed in practice, while QPR is reserved for applications demanding the most stringent conic‐error bounds regardless of runtime.\vspace{-2.5ex}
\subsection{Dynamic Performance}\label{res2}
We now evaluate the performance of the dynamic methods (DPR and DQPR) under a tight relaxation setting with $K_{\max}=5$. 
Table~\ref{tab:dynamic‐cuts} reports, for each instance, the average number of R\&F mappings \(\overline{\#\mathrm{R\&F}}\), the average number of outer cuts \(\overline{\#\mathrm{OC}}\) that actually take effect per SOC surface constraints, and the number of incumbent checks \(\#\mathrm{Check}\) performed. Note that \(\overline{\#\mathrm{OC}}\) counts only outer cuts in the first quadrant that remain active in the model, rather than the total cuts generated. Despite allowing up to $K_{\max}=5$ R\&F stages, most instances require fewer than one R\&F mapping per branch on average, and significantly fewer than the maximum number of outer cuts (\(2^{K_{\max}}+1\) cuts at even-split points in the first-quadrant arc). This confirms the effectiveness of lazy cut generation. DQPR avoids outer cuts entirely due to its retained SOC constraint, resulting in smaller models. Moreover, incumbent checks are only moderately frequent for both methods, further contributing to computational efficiency.
\begin{table}
  \centering\vspace{-3.5ex}
  \caption{Average Number of R\&F Mappings, Outer Cuts and Checks per SOC Surface Constraint at \(K_{\max}=5\)}
  \label{tab:dynamic‐cuts}
  \scriptsize
  \setlength{\tabcolsep}{4pt} 
  \begin{tabular}{%
      l   
      ccc|  
      cc  
    }
    \toprule
    \multirow{2}{*}{Instance} 
      & \multicolumn{3}{c}{DPR} 
      & \multicolumn{2}{c}{DQPR} \\[0.3ex]
    & \(\overline{\#\mathrm{R\&F}}\) 
      & \(\overline{\#\mathrm{OC}}\) 
      & \(\#\mathrm{Check}\)  
      & \(\overline{\#\mathrm{R\&F}}\) 
      & \(\#\mathrm{Check}\) \\
    \midrule
    IEEE 5  & 0 & 6.25 & 2 & 0.16 & 1\\
    IEEE 30& 0.47 & 4.45 & 16 & 0.61 & 9\\
    IEEE 118& 0.66 & 9.11 & 70 & 0.82 & 86\\
    IEEE 162 DTC& 1.00 & 10.48 & 77 & 1.14 & 64\\
    ACTIVSg200& 0.17 & 4.91 & 11 & 0.13 & 8\\
    IEEE 300& 0.39 & 6.99 & 35 & 0.45 & 34\\
    ACTIVSg500& 0.78 & 7.25 & 89 & 0.97 & 96\\
    GOC 793& 0.16 & 5.90& 37 & 0.12 & 31\\
    \bottomrule
  \end{tabular}\vspace{-3ex}
\end{table}
Tables~\ref{tab:pr_vs_dpr} and \ref{tab:qpr_vs_dqpr} compare PR versus DPR and QPR versus DQPR at \(K_{\max}=5\).  Column ``\(\Delta\)(\%)" reports the relative objective difference: \(\Delta = 100\times(\text{Obj}_{\text{dyn}} - \text{Obj}_{\text{stat}})/\text{Obj}_{\text{stat}},\) where \(\text{Obj}_{\mathrm{dyn}}\) and \(\text{Obj}_{\mathrm{stat}}\) denote the dynamic‐method and static‐formulation objectives, respectively. Across all eight test cases, DPR and DQPR closely match the objective values of their static counterparts within at most 0.06\% for DPR vs.\ PR, and 0.02\% for DQPR vs.\ QPR. These small gaps arise solely from the 0.1 \% optimality tolerance and confirm that the dynamic branch‐and‐cut procedure does not sacrifice optimality.
More striking are the solve‐time results. DPR achieves consistent speedups over PR in all but the smallest cases, with time reductions ranging from $1.2\times$ to $2.0\times$. This efficiency stems from selectively activating only the cuts needed for each branch. The QPR vs.\ DQPR comparison is more mixed: DQPR outpaces QPR on three of eight instances, but incurs longer solve times on the other five.  A likely explanation is that DQPR’s on‐demand cut generation incurs additional MISOCP node processing and incumbent checks, which can outweigh the benefit of a smaller model—especially when the static QPR already has few active cuts or the coarse relaxation repeatedly yields near‐feasible solutions, leading to higher dynamic overhead. 

\begin{table}
  \centering
  \caption{Comparison of Objective and Solve Time of PR vs. DPR at \(K_{\max}=5\)}
  \label{tab:pr_vs_dpr}
  \scriptsize
  \setlength{\tabcolsep}{4pt}
  \begin{tabular}{l
                  rrr|   
                  rrc   
                 }
    \toprule
    \multirow{2}{*}{Instance} 
      & \multicolumn{3}{c}{Obj} 
      & \multicolumn{3}{c}{Time (s)} \\
    & PR & DPR & $\Delta$ (\%) 
      & PR & DPR & Ratio \\
    \midrule
    IEEE 5  
      & 14\,999.69  
      & 14\,999.61  
      & $-0.00$  
      & \bf 0.04  
      & 0.07  
      & 0.58×  \\

    IEEE 30 
      & 6\,660.08  
      & 6\,660.07  
      & $-0.00$  
      & \bf 0.14  
      & 0.68  
      & 0.21×  \\

    IEEE 118 
      & 96\,392.43  
      & 96\,386.15  
      & $-0.01$  
      & 171.50  
      & \bf 84.80  
      & 2.02×  \\

    IEEE 162 
      & 101\,753.24  
      & 101\,697.09  
      & $-0.06$  
      & 1\,129.59  
      & \bf 757.87  
      & 1.49×  \\

    ACTIVSg 200 
      & 27\,470.17  
      & 27\,470.17  
      & $+0.00$  
      & 1.15  
      & \bf 0.72  
      & 1.59×  \\

    IEEE 300 
      & 550\,913.86  
      & 550\,850.34  
      & $-0.01$  
      & 1\,531.95  
      & \bf 994.55  
      & 1.54×  \\

    ACTIVSg 500 
      & 396\,818.06  
      & 396\,796.06  
      & $-0.01$  
      & 271.11  
      & \bf 221.20  
      & 1.23×  \\

    GOC 793 
      & 251\,311.09  
      & 251\,310.94  
      & $-0.00$  
      & 876.39  
      & \bf 476.31  
      & 1.84×  \\
    \bottomrule
  \end{tabular}\vspace{-3ex}
\end{table}
\begin{table}
  \centering\vspace{-3ex}
  \caption{Comparison of Objective and Solve Time of QPR vs. DQPR at \(K_{\max}=5\)}
  \label{tab:qpr_vs_dqpr}
  \scriptsize
  \setlength{\tabcolsep}{4pt}
  \begin{tabular}{l
                  rrr|   
                  rrc   
                 }
    \toprule
    \multirow{2}{*}{Instance} 
      & \multicolumn{3}{c}{Obj} 
      & \multicolumn{3}{c}{Time (s)} \\
    & QPR & DQPR & $\Delta$ (\%) 
      & QPR & DQPR & Ratio \\
    \midrule
    IEEE 5  
      & 14\,999.75  
      & 15\,001.41  
      & $+0.01$  
      & \bf 0.03  
      & 1.32  
      & 0.02×  \\

    IEEE 30 
      & 6\,662.23  
      & 6\,662.16  
      & $-0.00$  
      & \bf 0.11  
      & 3.56  
      & 0.03×  \\

    IEEE 118 
      & 96\,396.04  
      & 96\,397.21  
      & $+0.00$  
      & \bf 573.27  
      & 693.23  
      & 0.83× \\

    IEEE 162 
      & 101\,757.17  
      & 101\,732.23  
      & $-0.02$  
      & 1\,840.81  
      & \bf 1\,698.23  
      & 1.08× \\

    ACTIVSg 200 
      & 27\,470.31  
      & 27\,470.20  
      & $-0.00$  
      & \bf 4.79  
      & 10.20  
      & 0.47× \\

    IEEE 300 
      & 551\,688.23  
      & 551\,698.23  
      & $+0.00$  
      & 2\,206.00  
      & \bf 2\,012.20  
      & 1.10× \\

    ACTIVSg 500 
      & 396\,823.07  
      & 396\,843.06  
      & $+0.01$  
      & \bf 543.85  
      & 732.21  
      & 0.74× \\

    GOC 793
      & 251\,314.30  
      & 251\,309.30  
      & $-0.00$  
      & 1\,786.97  
      & \bf 1\,493.20  
      & 1.20× \\
    \bottomrule
  \end{tabular}\vspace{-3ex}
\end{table}\vspace{-3ex}
\subsection{Warm Start and Post-Processing}
To accelerate the MILP/MISOCP solves, we warm‐start formulation using a feasible solution from the IV formulation solved by IPOPT. Table \ref{tab:warm_start} reports, for each method on two representative cases (IEEE 118 and ACTIVSg 500) at \(K_{\text{max}}=5\), the time to find the first feasible solution $T_{\mathrm{FF}}$, the optimality gap of the first incumbent, the gap of the warm-start solution, and the total solve time both with and without warm start. 

Since PA is an approximation may lie outside its feasible region (Gap is marked “N/A”), and thus it cannot benefit from warm start—sometimes even experiencing slower solve times.  In contrast, relaxation‐based methods (PR, QPR, DPR, DQPR) accept the IV solution as a feasible incumbent, albeit with a modest initial gap, so they can begin pruning the search tree from a valid primal bound. In fact, obtaining a first incumbent can be quite time‐consuming (often dominates the total solution time), despite its typically good quality. By supplying the IV solution, typically computed within one second and with a small gap (often under 5 \%), we significantly accelerate speeds up the MILP/MISOCP solves for PR, QPR, and DPR. In contrast, DQPR sees little speedup (ratios near 1.0) due to the overhead of on-demand cut generation and incumbent checks.

After the branch‐and‐cut solve terminates, we apply a single LNS-based post‐processing step described in Section \ref{subsec:ald_det}. Table \ref{tab:lns_K2_detailed} shows, for each instance at \(K_{\max}=2\), pre-/post-processing conic error metrics and the post-processing time \(T_{\text{PP}}\). Whenever the fixed wedge still intersects the true cone, the LNS solutionconsistently improves both absolute and relative conic errors over the original DPR solution. If the wedge does not intersect the cone (e.g., IEEE 300, ACTIVSg 500), LNS is infeasible (marked “N/A”). Although not a relaxation and thus unable to provide a valid lower bound, the LNS procedure reliably recovers near-exact, cone-feasible solutions when a valid wedge exists. These solutions are effectively suboptimal QPR solutions, obtained with minimal overhead compared to full MISOCP solving.
\vspace{-2.5ex}
\subsection{Overall Performance}
\begin{figure*}\vspace{-3ex}
\centering
    \includegraphics[
    width=0.85\textwidth,
    height=2.5in
  ]{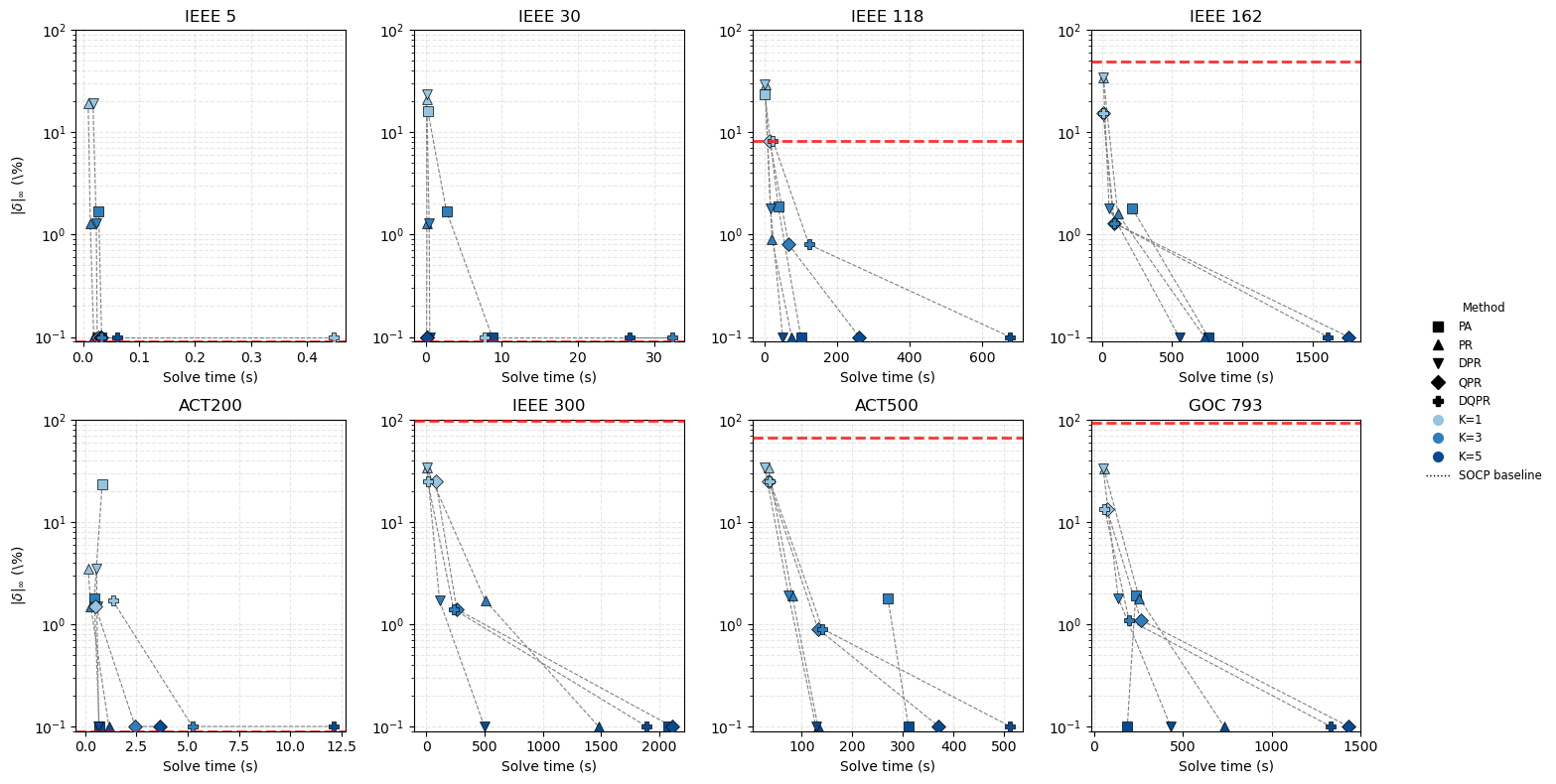}
    \caption{Solve‐time vs.\ maximum relative 4-D-cone violation $\|\bm \delta^{\text{rel}}\|_{\infty}$ (\%) for each method (K=1,3,5), across all eight test instances. Marker shape denotes method (square: PA; triangle up: PR; triangle down: DPR; diamond: QPR; pentagon: DQPR), while color denotes wedge depth K (light blue: 1; medium blue: 3; dark blue: 5). The red dashed line is the SOCP baseline error.}
    \label{fig:overall_scatter}\vspace{-3ex}
\end{figure*}

\begin{table}
  \centering
  \caption{Impact of Warm Start on First‐Feasible Time, Initial Gap, and Total Solve Time at \(K_{\max}=5\)}
  \label{tab:warm_start}
  \scriptsize
  \setlength{\tabcolsep}{3pt}
  \begin{tabular}{l l
                  r r r|   
                  r r r   
                 }
    \toprule
    \multirow{2}{*}{Method} 
      & \multirow{2}{*}{Instance} 
        & \multicolumn{3}{c}{No Warm Start} 
        & \multicolumn{2}{c}{With Warm Start}& \\
    &  & $T_{\mathrm{FF}}$ (s) & Gap\,(\%) & $T_{\mathrm{\text{total}}}$\, (s) 
      & Gap\,(\%) & $T_{\mathrm{\text{total}}}$\, (s) & Ratio \\
    \midrule
    PA  
      & IEEE 118  
      & 203     & 0.32   & \bf 214.8  
      & N/A     & 232.6  & 0.92\(\times\) \\
      & ACTIVSg 500
      & 281     & 0.01   & \bf 281.0  
      & N/A     & 311.4  & 0.90\(\times\) \\
    \midrule
    PR  
      & IEEE 118  
      & 171     & 0.10   & 171.5  
      & 0.78    &  \bf 75.6  & 2.26\(\times\) \\
      & ACTIVSg 500
      & 128     & 0.00   & 128.0  
      & 0.24    &  \bf 90.2  & 1.42\(\times\) \\
    \midrule
    QPR 
      & IEEE 118  
      & 352     & 0.23   & 573.2  
      & 0.77    & \bf 261.7  & 2.19\(\times\) \\
      & ACTIVSg 500
      & 543     & 0.00   & 543.9  
      & 2.32    & \bf 370.8  & 1.47\(\times\) \\
    \midrule
    DPR 
      & IEEE 118  
      &  72   & 0.37   &  84.8  
      & 2.91    &  \bf 50.9  & 1.66\(\times\) \\
      & ACTIVSg 500
      & 201   & 0.17   & 221.2  
      & 1.34    & \bf 189.2  & 1.17\(\times\) \\
    \midrule
    DQPR
      & IEEE 118  
      & 432   & 0.34   & 693.2  
      & 1.76    & \bf 677.1  & 1.02\(\times\) \\
      & ACTIVSg 500
      & 732     & 0.00   & \bf 732.2  
      & 2.58    & 745.2  & 0.98\(\times\) \\
    \bottomrule
  \end{tabular}\vspace{-3ex}
\end{table}
Figure~\ref{fig:overall_scatter} summarizes, for all eight benchmark instances, the trade‐off between solve time and maximum conic error \(\|\bm\delta\|_\infty\) for five primary methods (PA, PR, DPR, QPR, DQPR), each run with a warm start and three stage counts $K \in \{1,\,3,\,5\}$. Time is plotted linearly (x-axis), while conic error is log-scaled (y-axis). Marker shape indicates the method and color indicates the depth $K$. The SOCP baseline appears as a red dashed line.

Across all panels, SOCP solves in under 0.3 s but often produces very large violations (e.g. $\|\bm\delta\|_\infty>40\%$) on all but the smallest or easiest networks (IEEE 5, IEEE 30, ACT200).  PA often fails due to infeasibility (no point plotted). All relaxation-based methods remain feasible for every case and $K$. As $K$ increases, conic errors decrease steadily—by $K=5$, all methods reduce $\|\bm\delta\|_\infty$ below 0.1\%. Among proposed methods, PR and DPR show nearly identical accuracy, but DPR consistently solves faster, confirming the efficiency of dynamic cut generation. QPR and DQPR achieve the tightest relaxations but are slower due to MISOCP overhead, especially for $K=5$ the runtime gap widens substantially. 


These observations suggest the following practical recommendations: for small or easy networks, SOCP already attains sub-1 \% error within 1 s, making it the most efficient choice. For medium to large networks, when only moderate accuracy (e.g. $\|\bm \delta\|_\infty<5\%$) is required, DPR at $K=3$ is usually ideal. For example, on IEEE 118, DPR at $K=3$ runs in under 20 s while yielding $\|\bm\delta\|_\infty\approx 1.8\%$, whereas QPR or DQPR at $K=3$ take 60–140 s to reach similar or slightly lower errors. If one needs medium accuracy ($1\%>\|\bm\delta\|_\infty>0.1\%$), QPR at $K=3$ offer a good trade-offs: on IEEE 118, they achieve $\|\bm\delta\|_\infty\approx 0.8\%$ in roughly 67 s, which may be acceptable when tighter error is prioritized over speed. Finally, when high accuracy ($\|\bm\delta\|_\infty<0.1\%$) is required, all methods at $K=5$ meet this target; in that regime DPR is the fastest and therefore the best choice.
\begin{table}[t]
  \centering\vspace{-3.5ex}
  \caption{LNS Post‐Processing at \(K_{\max}=2\): Relative and Absolute Error Reduction, and Solve Time}
  \label{tab:lns_K2_detailed}
  \scriptsize
  \setlength{\tabcolsep}{4pt}
  \begin{tabular}{l
                  rr   
                  rr   
                  rr   
                  r    
                 }
    \toprule
    \multirow{2}{*}{Instance} 
      & \multicolumn{2}{c}{\(\|\bm\delta^{\text{rel}}\|_\infty\) (\%)} 
      & \multicolumn{2}{c}{\(\|\bm\delta^{\text{abs}}\|_\infty\)} 
      & \multicolumn{2}{c}{\(\|\bm\delta^{\text{abs}}\|_1\)} 
      & \multirow{2}{*}{\(T_{\text{PP}}\) (s)} \\
    & Before   & After    & Before    & After   & Before     & After   &  \\
    \midrule
    IEEE 5        & 4.81  & \(<0.1\)   & 0.67   & \(<0.1\)   & 1.85   & \(<0.1\)   & 0.12 \\
    IEEE 30       & 5.52  & 2.47       & 0.07   & 0.02       & 0.50   & 0.15       & 0.15 \\
    IEEE 118      & 6.80  & 3.33       & 2.48   & 1.21       & 12.16  & 3.92       & 0.48 \\
    IEEE 162      & 7.61  & 2.62       & 16.03  & 13.49      & 39.30  & 24.50      & 0.52 \\
    ACTIVSg 200   & 6.28  & 1.45       & 1.22   & 0.01       & 4.13   & 0.08       & 0.38 \\
    IEEE 300      & 7.50  & N/A        & 191.96 & N/A        & 644.08 & N/A        & N/A  \\
    ACTIVSg 500   & 7.91  & N/A        & 42.84  & N/A        & 237.12 & N/A        & N/A  \\
    GOC 793       & 7.91  & 4.87       & 97.68  & 0.38       & 171.46 & 0.85       & 0.55 \\
    \bottomrule
  \end{tabular}\vspace{-3ex}
\end{table}

In summary, SOCP works only on the smallest or easiest cases; PA is prone to infeasibility; PR offers a reasonable balance of speed and error; DPR consistently accelerates PR without loss of accuracy; and QPR/DQPR guarantee the tightest bounds at moderate cost. The scatter‐plot in Figure~\ref{fig:overall_scatter} clearly visualizes these trade‐offs across all eight test networks.

\vspace{-3ex}
\section{Conclusion}\label{sec:concl}
We proposed a unified framework for solving AC OPF globally. Central to this framework are two relaxations: PR and QPR, which approximate each branch-flow cone surface by a finite union of wedge-shaped regions and converge to exactness as the partition refines. To enable scalability, we developed dynamic branch-and-cut algorithms (DPR and DQPR) that incrementally tightens the SOC relaxation through on-the-fly cut generation. Warm starts and a lightweight LNS techniques further enhance solver performance and solution quality. Extensive experiments across eight PGLib–OPF benchmarks demonstrate that PR and DPR balance tightness and efficiency, QPR and DQPR provide higher accuracy when needed, and DPR remains the most practical choice for stringent accuracy requirements.

This framework opens several promising directions for future work. First, the dynamic cut-generation strategy could be extended to enforce nonconvex phase-angle recovery constraints which we neglect in this paper. Second, beyond feasibility, incorporating optimality-aware cut selection may further accelerate convergence. Third, exploiting structural dependencies among R\&F binary variables could inspire new cutting planes or enhance branching strategies. Finally, the framework readily generalizes to mixed-integer extensions of AC OPF, including optimal transmission switching
and AC unit commitment
, by embedding binary decisions within the same cut-generation architecture.

\appendices
\section{Proofs of Theorems and Propositions}
\subsection{Proof of Theorem \ref{thm:rnf}}
We omit the standard equivalence between the big-\(M\) reformulation and \(\lvert\cdot\rvert\) (see \cite{vielma2015mixed}) and focus on the R\&F recurrences. We prove both directions in turn.  

Let \((x,y,z)\in\Pi_K\).  By Definition~\ref{def:PiK}, on the plane \(z\) the point \((x,y)\) lies on some facet \(OV_nV_{n+1}\) of the inscribed \(2^{K+1}\)-pyramid, where
\[
V_n=\bigl(z\cos\varphi_n,\;z\sin\varphi_n\bigr),
\varphi_n=\tfrac{2(n-1)\pi}{2^{K+1}},\;n=1,\dots,2^{K+1}.
\]
Hence \((x,y)\) is a convex combination of \(V_n\) and \(V_{n+1}\).
By construction of the R\&F mapping, setting 
\((g_0,h_0)=(|x|,|y|)\) and then iteratively folds and rotates that same convex combination into the first quadrant and into the final small sector.  Because \((x,y)\) lay on the facet \(V_nV_{n+1}\) in \(z\)-plane, one shows by direct trigonometry that
\[
g_K=z\cos\theta_K,
\qquad
0\le h_K\le z\sin\theta_K,
\]
which is exactly \eqref{rffinal}. Thus \((x,y,z)\) satisfies \eqref{rf1}, \eqref{rfmap}, and \eqref{rffinal}.
 
Conversely, suppose \((x,y,z)\) satisfies the R\&F system \(\eqref{rf1}\)–\(\eqref{rffinal}\).  Then by \(\eqref{rffinal}\),
\[
g_K = z\cos\theta_K,
\qquad
0\le h_K\le z\sin\theta_K,
\]
so the point \((g_K,h_K)\) lies on the line segment between 
\(\tilde{V}_1=(z\cos\theta_K,0)\) 
and 
\(\tilde{V}_2=(z\cos\theta_K,z\sin\theta_K)\).  Hence there is \(\lambda\in[0,1]\) with
\[
\begin{pmatrix}g_K\\h_K\end{pmatrix}
=(1-\lambda)\tilde{V}_1+\lambda\tilde{V}_2.
\]
Now invert the \(K\) R\&F steps and the initial quadrant fold (using the known sign–binaries from the big-\(M\) linearization) to recover \((x,y)\) from \((g_K,h_K)\).  This shows
\((x,y)\) is a (signed) convex combination of exactly two adjacent vertices of the base pyramid,
\[
(x,y)
=(1-\lambda)V_n + \lambda\,V_{n+1},
\]
so \((x,y,z)\) satisfies the SOS2 representation \eqref{eq:PA-weights}–\eqref{vceq}.  By Definition~\ref{def:PiK}, \((x,y,z)\in\Pi_K\).

Combining both directions completes the proof of Theorem~\ref{thm:rnf}.\IEEEQEDhere
\subsection{Proof of Theorem \ref{thm:eps}}
We prove each claim by showing that any feasible point \((x,y,z)\) of the given \(K\)-stage formulation satisfies \[
-\varepsilon \;\le\;
\frac{x^2+y^2 - z^2}{z^2}
\;\le\;\varepsilon
\]
with the stated \(\varepsilon\). First observe that axial‐symmetry and R\&F folds preserve the Euclidean norm of \((x,y)\), so 
\[
g_K^2 + h_K^2 \;=\; x^2 + y^2.
\]
\paragraph{\(K\)-PA is \(\sin^2\theta_K\)-feasible} In \(K\)-PA the final facet \eqref{rffinal} imposes
\[
g_K 
= z\cos\theta_K,
\quad
0\le h_K\le z\sin\theta_K,
\]
so
\[
z^2\cos^2\theta_K \le g_K^2+h_K^2 
\le z^2\cos^2\theta_K+z^2\sin^2\theta_K = z^2.
\]
Thus
\[
-\sin^2\theta_K 
\le 
\frac{x^2+y^2 -z^2}{z^2} 
\le 0
\]
so the relative conic error \(\delta=|\frac{x^2+y^2 -z^2}{z^2}|\le \sin^2\theta_K\).
\paragraph{\(K\)-PR is \(\tan^2\theta_{K+1}\)-feasible} Recall the three boundaries of the final relaxations in the \(z\)‐plane \eqref{eqpr}.
By Cauchy–Schwarz and the third cut,
\[\begin{aligned}
    g_K^2+h_K^2&=(g_K^2+h_K^2)(\cos^2\theta_{K+1}+\sin^2\theta_{K+1})\\
    &\ge(g_K\cos\theta_{K+1} + h_K\sin\theta_{K+1})^2\\
 &\ge
 z^2\cos^2\theta_{K+1},
\end{aligned}
\] hence
\[
 \frac{x^2+y^2 -z^2}{z^2} \ge -\sin^2\theta_{K+1}.
\]
To obtain the upper bound on \(x^2+y^2=g_K^2+h_K^2\), consider the quadratic program 
\[
\begin{aligned}
\max_{g_K,h_K\ge0}\quad & g_K^2 + h_K^2\\
\text{s.t.}\quad 
& g_K - z \;\le\;0,\\
& g_K\cos\theta_K + h_K\sin\theta_K - z\;\le\;0,\\
& -\,g_K\cos\theta_{K+1}-h_K\sin\theta_{K+1} + z\cos\theta_{K+1}\;\le\;0.
\end{aligned}
\]
Let \(R\) denote its feasible region.  Since \(R\) is a nonempty compact convex set and \(f(g,h)=g^2+h^2\) is a continuous convex function, Bauer’s maximum principle implies that the maximum of \(f\) over \(R\) occurs at one of \(R\)’s extreme points. Direct enumeration shows these vertices are
\[
\bigl(z,0\bigr),\quad
\bigl(z\cos\theta_K,\;z\sin\theta_K\bigr),\quad
\bigl(z,\;z\tan\theta_{K+1}\bigr).
\]
Evaluating \(f\) at each vertex shows that
\[
\max_{(g,h)\in R}(g^2+h^2)=z^2 + z^2\tan^2\theta_{K+1}
= z^2\sec^2\theta_{K+1},
\]
and therefore
\[
\frac{x^2+y^2-z^2}{z^2}
\le
\sec^2\theta_{K+1}-1=
\tan^2\theta_{K+1}.
\]
Combining these bounds shows 
\[-\sin^2\theta_{K+1}\le \frac{x^2+y^2-z^2}{z^2}\le\tan^2\theta_{K+1},\]
hence \(\delta=|\frac{x^2+y^2 -z^2}{z^2}|\le \tan^2\theta_{K+1}\).
\paragraph{\(K\)-QPR is \(\sin^2\theta_{K+1}\)-feasible.}  
In \(K\)-QPR we retain the exact SOC constraint \(x^2+y^2\le z^2\) and add only the inner‐cut
\[
g_K\cos\theta_{K+1} + h_K\sin\theta_{K+1} \ge z\cos\theta_{K+1}.
\]
Thus
\[
x^2+y^2 \;=\; g_K^2+h_K^2 \;\le\; z^2
\;\Longrightarrow\;
\frac{x^2+y^2-z^2}{z^2}\;\le\;0,
\]
and by the same Cauchy–Schwarz argument as above,
\[
\frac{x^2+y^2 - z^2}{z^2}
\;\ge\;
-\sin^2\theta_{K+1}.
\]
Hence \(\lvert\delta\rvert\le\sin^2\theta_{K+1}\), completing the proof.\IEEEQEDhere

\subsection{Proof of Proposition~\ref{prop:consist}}
Fix \(K\ge0\) and consider the static \(K\)-PR feasible region \(\mathcal X_K\) and the final feasible region at termination of the \(K\)‐DPR algorithm \(\mathcal X_{\rm dyn}\). By construction, every cut that DPR ever introduces is one of the static \(k\)-PR cuts for some $k\le K$, and since
\(\mathcal X_K\subseteq\mathcal X_{K-1}\subseteq\cdots\subseteq\mathcal X_0,\)
we have
\[
\mathcal X_K \;\subseteq\;\mathcal X_{\rm dyn}
\]
throughout the run.  Upon termination, the incumbent \(x^*\) violates no further static $K$-PR cuts, so \(x^*\in\mathcal X_K\).

Because DPR minimizes the objective over \(\mathcal X_{\rm dyn}\), its final incumbent satisfies
\[
f(\mathbf{x}^*) \;=\;\min_{\mathbf{x}\in \mathcal X_{\rm dyn}} f(\mathbf{x})
\;\le\;
\min_{\mathbf{x}\in \mathcal X_K} f(\mathbf{x}).
\]
On the other hand, \(\mathbf{x}^*\in\mathcal X_K\) implies
\(\min_{\mathbf{x}\in\mathcal X_K}f(\mathbf{x})\le f(\mathbf{x}^*)\).
Hence equality holds, and \(\mathbf{x}^*\) also attains the minimum over \(\mathcal X_K\).  That is, \(\mathbf{x}^*\) is optimal for the static $K$-PR problem. An identical argument shows $K$-DQPR returns the static $K$-QPR optimum.

\subsection{Proof of Proposition~\ref{thm:bound}}
We prove by induction on \(n\) that for each \(n=0,1,\dots,K\),
\[
g_n\cos\theta_n + h_n\sin\theta_n \le z,
\]
given the base bounds
\[
0 \;\le\; g_0,h_0 \;\le\; z,
\]
and the R\&F recurrences
\[
\begin{aligned}
g_n &= \cos\theta_n\,g_{n-1} +\sin\theta_n\,h_{n-1},\\
h_n &= \bigl|-\sin\theta_n\,g_{n-1} + \cos\theta_n\,h_{n-1}\bigr|,\\
g_n &\le z,
\end{aligned}
\quad n=1,2,\dots,K.
\]
\paragraph{Base case \(n=0\)}
Here \(\theta_0=\pi/2\), so
\[
g_0\cos\theta_0 + h_0\sin\theta_0
= g_0\cdot0 + h_0\cdot1
= h_0
\le z,
\]
by the assumption \(h_0\le z\).
\paragraph{Inductive step}
Assume for some \(n\ge0\) that
\[
g_n\cos\theta_n + h_n\sin\theta_n \le z.
\]
We must show
\(\;g_{n+1}\cos\theta_{n+1} + h_{n+1}\sin\theta_{n+1}\le z.\)

From the definitions,
\[
\begin{aligned}
&g_{n+1}\cos\theta_{n+1} + h_{n+1}\sin\theta_{n+1}\\
=&
\bigl(\cos\theta_{n+1}\,g_n + \sin\theta_{n+1}\,h_n\bigr)\cos\theta_{n+1}\\
+&\bigl|-\sin\theta_{n+1}\,g_n + \cos\theta_{n+1}\,h_n\bigr|\sin\theta_{n+1}.
\end{aligned}
\]
We consider two cases:
1. \emph{Inner argument nonnegative:}
   \[
   -\sin\theta_{n+1}\,g_n + \cos\theta_{n+1}\,h_n
   \ge 0.
   \]
   Then
   \[
   \begin{aligned}
   &g_{n+1}\cos\theta_{n+1} + h_{n+1}\sin\theta_{n+1}\\
   =& g_n(\cos^2\theta_{n+1}-\sin^2\theta_{n+1})
     +2h_n\sin\theta_{n+1}\cos\theta_{n+1}\\
   =& g_n\cos(2\theta_{n+1}) + h_n\sin(2\theta_{n+1})\\
   =& g_n\cos\theta_n + h_n\sin\theta_n
   \;\le\; z,
   \end{aligned}
   \]
   using \(2\theta_{n+1}=\theta_n\) and the induction hypothesis.

2. \emph{Inner argument negative:}
   \[
   -\sin\theta_{n+1}\,g_n + \cos\theta_{n+1}\,h_n
   <0.
   \]
   Then \(\lvert\cdot\rvert\) flips sign and
   \[
   \begin{aligned}
   &g_{n+1}\cos\theta_{n+1} + h_{n+1}\sin\theta_{n+1}\\
   =& g_n(\cos^2\theta_{n+1}+\sin^2\theta_{n+1})
   = g_n
   \le z,
   \end{aligned}
   \]
   since \(g_n\le z\) is enforced at every stage.  
   
In both cases we conclude
\[g_{n+1}\cos\theta_{n+1} + h_{n+1}\sin\theta_{n+1}\le z,\]
completing the induction.  Setting \(n=K\) yields the desired
\(\;g_K\cos\theta_K + h_K\sin\theta_K\le z.\)  \IEEEQEDhere

\section{Dynamic Relaxation Algorithm Framework}\label{appendix-alg} $\mathbf{k}$ is a $2\times L$ matrix that stores the number of current R\&F depth, where the rows correspond to the two types of conical surface constraints and the columns correspond to branches. $\mathbf{x}^*=(x_1^*,...,x_L^*)$ is an incumbent, where $x_l^*$ is the component of $x^*$ regarding branch $l$. $x_l^*=(\hat{x},\hat{y},\hat{z},\hat{\bm \beta}_0,\hat{g}_0,\hat{h}_0,...,\hat{\beta}_k,\hat{g}_k,\hat{h}_k,0,..,0)$ omits the subscripts of constraint types and branches for simplicity.

  \begin{breakablealgorithm}
   \caption{{\bf Dynamic Relaxation Framework} for AC OPF}
   \setlength{\algorithmicindent}{0.5em}
   \begin{algorithmic}[1]\label{alg-1}
   \REQUIRE AC OPF instance \(\mathcal P\), 
    initial depth \(K_{\mathrm{init}}\), 
    maximum depth \(K_{\max}\).
     \ENSURE 
    Stored incumbents $\mathcal F$, best MILP solution $\mathbf{x}^{\mathrm{opt}}$, optional SOCP refined solution $\mathbf{x}^{\mathrm{LNS}}$.
    \STATE  \(\mathcal C \leftarrow\) all constraints of the \(K_{\mathrm{init}}\)-stage relaxation (PR or QPR)
    \STATE \(\mathcal F \leftarrow \varnothing\)
\FOR{cone‐type \(i=1,2\) and branch \(l=1,\dots,L\)}
\STATE \(k_{i,l}\leftarrow K_{\mathrm{init}}\)
\ENDFOR
\STATE \texttt{--- Optional Warm Start ---}
\STATE  \({\mathbf{x}^{\mathrm{warm}}}\leftarrow\) solve the IV‐formulation by IPM
\STATE \(\mathcal F\leftarrow\mathcal F\cup\{\mathbf{x}^{\mathrm{warm}}\}\)
\STATE \textbf{Define} $\mathcal P'$ = problem with objective from $\mathcal P$ and constraints set $\mathcal C$.
  \STATE \textbf{Begin branch‐and‐cut on} $\mathcal P'$:
\WHILE{solver not terminated}
    \IF{new integer‐feasible incumbent $\mathbf{x}^*$ is found}
    \STATE \(flag \leftarrow \)\texttt{true}
            \FOR{$i=1:2,l = 1:L$}
                \IF{\Call{ViolationCheck}{\(x^*_l,i,l,\mathbf k\)}}
                \IF{$x_l^*$ lies inside cone}
                    \STATE \Call{AddInnerCut}{\(x_l^*, i,l,\mathbf k,\mathcal C\)}; \(flag\gets\texttt{false}\)
                \ELSIF{$x_l^*$ lies outside cone}
                    \STATE \Call{AddOuterCut}{\(x_l^*, i,l,\mathbf k,\mathcal C\)}; \(flag\gets\texttt{false}\)
                \ENDIF
                \ENDIF
            \ENDFOR
            \IF{flag}
                \STATE $\mathcal{F}\leftarrow\mathcal{F}\cup\{\mathbf{x}^*\}$
            \ENDIF
    \ENDIF
\ENDWHILE
\STATE \(x^{\mathrm{opt}}\leftarrow\) best incumbent found
\STATE \texttt{--- Optional post‐processing (LNS) ---}
\FOR{$i=1:2,l = 1:L$}
     \STATE Fix all R\&F binaries \(\bm\beta_{0},\ldots, \beta_{K_{\max}}\) to values in \(\mathbf{x}^{\mathrm{opt}}\)
    \STATE $\mathcal{C} \leftarrow\mathcal{C} \cup \{\text{SOC constraints for } i,l\}$
\ENDFOR
 \STATE Resolve as pure SOCP to obtain refined solution \(\mathbf{x}^{\mathrm{LNS}}\)
\RETURN \(\mathcal F,\mathbf{x}^{\mathrm{opt}},\mathbf{x}^{\mathrm{LNS}}\)
   \end{algorithmic}
  \end{breakablealgorithm}

\begin{algorithm}
  \caption{\textsc{ViolationCheck} subroutine}
  \label{alg-NeedsInnerCut}
\begin{algorithmic}[1]
  \REQUIRE incumbent element \(x_l^*\), cone-type \(i\), branch $l$, depths \(\mathbf{k}\)
  \ENSURE \texttt{true} if a cut is required
  \STATE extract and record \(\hat{\beta}_k,\hat{g}_{k},\hat{h}_{k}\) for \(k=k_{i,l}+1,\ldots,K_{\max}\) via R\&F \eqref{rfmap} from \(x^*_l\)
  \IF{\((\hat{g}_{K_{\max}},\hat{h}_{K_{\max}},\hat z)\) violate cuts \eqref{eqpr}}
    \RETURN \texttt{true}
  \ELSE
    \RETURN \texttt{false}
  \ENDIF
\end{algorithmic}
\end{algorithm}

\begin{algorithm}
  \caption{\textsc{AddInnerCut} subroutine}
  \label{alg-AddInnerCut}
\begin{algorithmic}[1]
  \REQUIRE incumbent element \(x_l^*\), cone‐type \(i\), branch \(l\), depths \(\mathbf k\), cut‐set \(\mathcal C\)
  \ENSURE append R\&F mappings and inner cut to \(\mathcal C\)
  \STATE \(k_{\rm old}\gets k_{i,l}\)
  \STATE $k_\textrm{new}\leftarrow\min\{k:\hat{g}_k,\hat{h}_k,\hat z \text{ violate inner cut of \eqref{eqpr}}\}$.
  \FOR{\(k = k_{\rm old}+1:k_{\rm new}\)}
    \STATE \(\mathcal C \leftarrow\mathcal C \cup \{\text{R\&F mapping \eqref{rfmap} at level } k \text{ for } i,l\}\)
  \ENDFOR
  \STATE \(\mathcal C \leftarrow\mathcal C \cup \{k\text{-inner cut of \eqref{eqpr} for } i,l\}\)
  \STATE \(k_{i,l}\gets k_{\rm new}\)
\end{algorithmic}
\end{algorithm}

\begin{algorithm}
  \caption{\textsc{AddOuterCut} subroutine}
  \label{alg-AddOuterCut}
\begin{algorithmic}[1]
  \REQUIRE incumbent component \(x_l^*\), cone‐type \(i\), branch \(l\), depths \(\mathbf k\), cut‐set \(\mathcal C\)
  \ENSURE append the next outer‐cut(s) to \(\mathcal C\)
  \STATE Compute $\delta_l(x^*_l)$ and $\arg(\hat{g}_{k_{i,l}},\hat{h}_{k_{i,l}},\hat z)$.
  \STATE $k_{\text{new}}\leftarrow\min\{k:\delta_l(x^*_l)>\delta_k^{out}\}$.
    \STATE $\psi_1\leftarrow\frac{\pi}{2^{k_{\text{new}}+1}}\lfloor\arg(\hat{g}_{k_{i,l}},\hat{g}_{k_{i,l}},z)/(\frac{\pi}{2^{k_{\text{new}}+1}})\rfloor$,\\
    $\psi_2\leftarrow\frac{\pi}{2^{k_{\text{new}}+1}}\lceil\arg(\hat{g}_{k_{i,l}},\hat{g}_{k_{i,l}},z)/(\frac{\pi}{2^{k_{\text{new}}+1}})\rceil$
    \STATE $\mathcal{C}\leftarrow\mathcal{C}\cup \{\text{outer cuts \eqref{otc} of }\psi_1 $ and $ \psi_2\}.$
\end{algorithmic}
\end{algorithm}

\vspace{-2ex}
\section*{Acknowledgment}
This research was supported by the National Key R\&D Program of China (No.
2022YFB2403400) and the Chinese NSF grants (No. 12201620). 
The computations were done on the high performance computers of State Key Laboratory of Mathematical Sciences.
\vspace{-2ex}
\bibliographystyle{IEEEtran}
\bibliography{bibtex/bare_jrnl_new_sample4.bib}{}

\newpage

 




\vfill

\end{document}